\documentclass[11pt,letterpaper]{amsart}
\usepackage[margin=2.5cm]{geometry}
\usepackage{setspace}
\usepackage{amsmath}
\usepackage{mathtools}
\usepackage{amssymb}
\usepackage{amsthm}
\usepackage{todonotes}
\usepackage{bbm}
\usepackage{stmaryrd}
\usepackage[normalem]{ulem}
\usepackage{mathdots}
\usepackage{float}
\usepackage[all]{xy}

\usepackage[pagebackref,colorlinks=true,linkcolor=cyan,citecolor=magenta]{hyperref}

\usepackage[alphabetic,nobysame, initials, backrefs]{amsrefs}
\usepackage{tikz-cd}
\usepackage{tikz}
\usepackage{nicematrix}
\usepackage{xcolor}
\newcommand{\rev}[1]{{\color{red}#1}}

\tikzset{
  symbol/.style={
    draw=none,
    every to/.append style={
      edge node={node [sloped, allow upside down, auto=false]{$#1$}}}
  }
} 

\usepackage{imakeidx}
\makeindex[columns=2]

\newcommand{\IC}{\mathcal{IC}}

\newcommand{\Z}{\mathbb{Z}}

\newcommand{\BC}{\mathbb{C}}

\newcommand{\ol}{\overline}

\newcommand{\SL}{\mathrm{SL}}
\newcommand{\GL}{\mathrm{GL}}

\newcommand{\dualgroup}[1]{{#1}^\vee}
\newcommand{\irrepdual}[1]{\widehat{#1}}
\newcommand{\Lgroup}[1]{{\hskip-2 pt \,^L\hskip-1pt{#1}}}

\newcommand{\pure}{\mathrm{pure}}

\newcommand{\pr}{\mathrm{pr}}

\newcommand{\mfr}[1]{\mathfrak{#1}}

\newcommand{\Fr}{\mathrm{Fr}}

\newcommand{\Ind}{\mathrm{Ind}}

\newcommand{\LLC}{\mathrm{LLC}}
\newcommand{\Per}{\mathrm{Per}}

\newcommand{\temp}{\mathrm{temp}}

\newcommand{\Lift}{\mathrm{Lift}}

\newcommand{\ABV}{{\mbox{\raisebox{0pt}{\scalebox{0.6}{$\mathrm{ABV}$}}}}}

\newcommand{\locsys}[2]{\underline{#1}_{#2}}
\newcommand{\stsh}[2]{{#1}_{#2}}

\newcommand{\Irrep}{\mathrm{Irrep}}

\renewcommand{\implies}{\Rightarrow}

\newcommand{\comment}[1]{}

\newcommand{\CC}{\mathbb{C}}

\newtheorem{thm}{Theorem}[section]

\newtheorem{lemma}[thm]{Lemma}
\newtheorem{prop}[thm]{Proposition}

\newtheorem {ques/conj}[thm]{Question/Conjecture}
\newtheorem {hypothesis}[thm]{Hypothesis}

\newtheorem{defn}[thm]{Definition}
\newtheorem{remark}[thm]{Remark}

\newtheorem{desideratum}[thm]{Desideratum}

\newtheorem*{globalconj*}{Global Conjecture}
\newtheorem*{localconj*}{Local Conjecture}

\newtheorem*{nonzero*}{Conjecture on the non-vanishing of the normalized intertwining operators}
\newtheorem*{holo*}{Conjecture on the holomorphicity of the normalized intertwining operators}

\DeclareMathOperator{\Ad}{Ad}

\DeclareMathOperator{\Hom}{Hom}
\DeclareMathOperator{\Gal}{Gal}

\let\oldbullet\bullet
\renewcommand{\bullet}{{\vcenter{\hbox{\tiny$\oldbullet$}}}}

\numberwithin{equation}{section}

\begin{document}

\title[induction restriction]{Parabolic induction and pullback of perverse sheaves on vogan varieties}

\author[Lo]{Chi-Heng Lo}
\address{Department of Mathematics,
National University of Singapore,
119076, Singapore}
\email{ch\_lo@nus.edu.sg}

\author[Ray]{Mishty Ray}
\address{Department of Mathematics\\
University of British Columbia\\
Vancouver, BC, V6T 1Z4,
CANADA}
\email{mishtyray@math.ubc.ca}


\date{\today}


\thanks{}

\begin{abstract}
Let $M$ be a Levi subgroup of a quasisplit connected reductive group $G$ over a $p$-adic field. Vogan's geometric perspective on the local Langlands correspondence fixes a perfect pairing between the Grothendieck group of representations of $G$ and its pure inner forms with a fixed infinitesimal parameter and the Grothendieck group of equivariant perverse sheaves on the associated Vogan variety. Granting a list of desiderata for the enhanced local Langlands correspondence, we show that parabolic induction $\Ind_M^G$ is the transpose of geometric restriction $\xi^*$ with respect to this pairing if and only if the $p$-adic analogue of the Kazhdan--Lusztig hypothesis and the expected multiplicity formula for tempered induction both hold. The identity refines to one for each pure inner form of $G$. \end{abstract}


\maketitle

\tableofcontents

\section{Introduction}\label{sec: introduction}

Let $F$ be a nonarchimedean local field of characteristic $0$ and let $G$ be a quasisplit connected reductive group over $F$.  The local Langlands correspondence predicts a finite-to-one parametrization of the smooth irreducible representations of $G(F)$ by $L$-parameters $\phi: W_F\times\SL_2(\BC)\to\Lgroup{G}$ where $\Lgroup{G}=\dualgroup{G}\rtimes W_F$ and $\dualgroup{G}$ is the $\CC$-points of the dual group of $G$.  An enhancement of the local Langlands correspondence refines the finite-to-one map into a bijection using irreducible representations of a component group attached to $\phi$. More precisely, on one side, we consider the set $\Pi^{\pure}(G)$ of equivalence classes of smooth irreducible representations of $G(F)$ along with those of its pure inner forms $G_\delta(F)$, as $\delta$ varies over elements of the Galois cohomology group $H^1(F,G)$. On the other hand, we have $\Phi^{\mfr{e}}(G)$, the set of equivalence classes of pairs $(\phi,\tau)$, where $\phi$ is a Langlands parameter and $\tau$ is an irreducible representation of the component group $A_\phi:=\pi_0(Z_{\dualgroup{G}}(\phi))$. This formulation is due to Vogan; see \cite{Vog93}, or for the archimedean case, \cite{ABV92}. We denote the enhanced local Langlands correspondence by
\[ \LLC^{\mfr{e}}: \Pi^{\pure}(G) \to \Phi^{\mfr{e}}(G).\]
The enhanced local Langlands correspondence is constructed for $\GL_n(F)$ in \cites{HT01, Hen00, scholze}, for classical groups in \cites{Art13, mokendoscopic2015, KMSW14, MR18, Ish24}, and for $G_2$ in \cites{GanSavin2023, aubert2023explicitlocallanglandscorrespondence}. There has been far more progress towards this conjecture than we are able to cite here; we move on for the sake of brevity.
A central aspect of this program is functoriality. An $L$-homomorphism $\Lgroup{H}\to\Lgroup{G}$ should give rise to a transfer of (packets of) representations from $H$ to $G$.  In \cite{ABV92}*{Proposition 1.37}, the authors formulate a Langlands functoriality map for archimedean groups in a geometric setting and utilize it when $H$ is an endoscopic group of $G$. This is extended to $p$-adic groups in \cite{Vog93} and reformulated in the language of Vogan varieties in \cite{CFMMX22}. We utilize the latter formulation in reviewing the geometric setting below. 

Given a Langlands parameter $\phi$, we attach an infinitesimal parameter by
\[ \lambda_{\phi}:= \phi\left(w,\begin{pmatrix}
    |w|^{1/2}&\\&|w|^{-1/2}
\end{pmatrix}\right).\]
The map $\phi \mapsto\lambda_{\phi}$ is a finite-to-one surjection modulo equivalence.  Let $\Phi_{\lambda}(G)$ be the fibre of this surjection at $\lambda$. We set $\Phi^{\mathfrak{e}}_\lambda(G)$ to be the subset of $\Phi^{\mfr{e}}(G)$ consisting of $(\phi,\tau)$ such that $\lambda_\phi=\lambda$ (or its $\dualgroup{G}$-conjugate). This enables us to start building the geometric spaces we need for Vogan's perspective.
Fixing an infinitesimal parameter $\lambda$ of $G$, we obtain the prehomogeneous vector space $V_\lambda$ called the {Vogan variety} with an action of the group $H_\lambda$ (see Section \ref{subsec: a moduli space for langlands parameters} for definitions).  Every $L$-parameter $\phi$ with infinitesimal parameter $\lambda$ determines an $H_\lambda$-orbit $C_\phi$, the representation $\tau$ determines an equivariant local system $\locsys{\tau}{C_\phi}$ on $C_\phi$, and $\mathcal{IC}(\locsys{\tau}{C_\phi})$ is the corresponding simple object of the category $\Per_{H_\lambda}(V_\lambda)$ of $H_\lambda$-equivariant perverse sheaves on $V_\lambda$. The enhanced local Langlands correspondence gives the bijections:
  \[\begin{array}{ccccc}
\Pi_\lambda^\pure(G) &{\longrightarrow} &\Phi^{\mathfrak{e}}_\lambda(G) & {\longrightarrow} &  \Per_{H_\lambda}(V_\lambda)^{\mathrm{simple}}_{\mathrm{/iso}} \\
\pi & {\longmapsto} & (\phi, \tau) & {\longmapsto} & \mathcal{IC}(\locsys{\tau}{C_\phi})=: \mathcal{P}(\pi).
\end{array}
\]
The above bijection defines a perfect pairing of Grothendieck groups
\[
    K\Pi^\pure_\lambda(G) \times K\Per_{H_\lambda}(V_\lambda) \to \mathbb{Z},
\]
\begin{equation}\label{equation: pairing in the introduction}
\langle \pi,\mathcal{P}\rangle =
\begin{cases}
    (-1)^{d(\pi)}e(\pi) & \text{if } \mathcal{P} = \mathcal{P}(\pi),\\
    0 & \text{otherwise},
\end{cases}
\end{equation}
where $d(\pi) \coloneqq \dim C_{\phi}$ and $e(\pi)$ is the Kottwitz sign of the pure inner form associated to $\pi$; see \eqref{equation: kottwitz sign}.

Returning to the question of functoriality, in this paper we only consider $H=M$, a Levi subgroup of $G$. Let $\xi: {}^L M \hookrightarrow {}^L G$ be the $L$-embedding. Fix an infinitesimal parameter $\lambda_M$ of $M$ and put $\lambda=\xi\circ\lambda_M$. We have the  $\Z$-linear map
\[\Ind_M^G:K\Pi^\pure_{\lambda_M}(M)\to K\Pi^\pure_{\lambda}(G),\qquad \pi \mapsto \Ind_P^G(\pi),\]
where $P$ is the standard parabolic subgroup of $G$ with Levi component $M$ and $\Ind_P^G$ is normalized parabolic induction (for a representation of a pure inner form of $M$ one induces to the corresponding pure inner form of $G$; see Section \ref{subsec: parabolic induction of reps linear map}). The map $\xi$ induces an inclusion $H_{\lambda_M} \hookrightarrow H_\lambda$ and an equivariant inclusion $V_{\lambda_M} \hookrightarrow V_\lambda$. We obtain the restriction functor $\xi^*:D_{H_\lambda}(V_\lambda) \to D_{H_{\lambda_M}}(V_{\lambda_M})$ between the equivariant bounded derived categories of constructible sheaves, and thus a $\Z$-linear map
\[\xi^*: K\Per_{H_\lambda}(V_\lambda) \to K\Per_{H_{\lambda_M}}(V_{\lambda_M})\]
between the respective Grothendieck groups. Then we define ${}^t\xi^* : K\Pi^\pure_{\lambda_M}(M)\to K\Pi^\pure_{\lambda}(G)$ by
\[ \langle {}^t\xi^*(\pi), \mathcal{P}\rangle=\langle \pi, \xi^*(\mathcal{P})\rangle\]
for any $\pi\in K\Pi_{\lambda_M}^\pure(M)$ and any $\mathcal{P}\in K\Per_{H_\lambda}(V_\lambda)$.
In this paper we establish the equality
 \begin{align}\label{eq ind res intro}
     \Ind_M^G={}^t\xi^*,
 \end{align}
 assuming specific desiderata from the local Langlands correspondence. Below, we expand on two of them, since they are of special importance to a converse statement. 

The first concerns the decomposition of tempered induction. If $(\phi_M, \tau_M)$ is an enhanced $L$-parameter of $M$ with $\phi_M$ tempered, and $\phi=\xi\circ\phi_M$, then one expects
\begin{align}\label{eq temp decomp intro}
\Ind_P^G\big(\pi(\phi_M,\tau_M)\big)=\bigoplus_{\tau\in \irrepdual{A}_\phi} m(\tau_M,\tau)\,\pi(\phi,\tau),
\end{align}
where $m(\tau_M,\tau)$ is the multiplicity of $\tau_M$ in the restriction of $\tau$ to $A_{\phi_M}$. This is Desideratum \ref{desideratum: induction of tempered representations}. 

The second is the $p$-adic analogue to the Kazhdan-Lusztig hypothesis. The Grothendieck group $K\Pi_\lambda^\pure(G)$ (resp. $K\Per_{H_\lambda}(V_\lambda)$) has two bases, one consists of irreducible representations $\pi(\phi,\tau)$ (resp. simple perverse sheaves $\mathcal{IC}(\locsys{\tau}{C_\phi})$), and the other consists of standard modules $\Delta(\pi)$ of $\pi$ (resp. standard sheaves $\stsh{\tau}{C_{\phi}}$). We let $m_{\lambda}$ and $c_{\lambda}$ be the change of basis matrices between the simple and standard objects in $K\Pi_\lambda^\pure(G)$ and $K\Per_{H_\lambda}(V_\lambda)$, respectively. The $p$-adic Kazhdan-Lusztig hypothesis states that 
\[ m_{\lambda}= {}^tc_{\lambda}.\]
Equivalently, the evaluation of the pairing \eqref{equation: pairing in the introduction} on standard objects is given by 
\begin{equation}\label{equation: pairing on standards in the introduction}
    \begin{array}{rcl}
        \langle \Delta(\pi), \stsh{\tau}{C_{\phi}} \rangle 
            &\mapsto& 
            \begin{cases}
                e(\pi) & \LLC^{\mathfrak{e}}_G(\pi) = (\phi,\tau), \\
                0 & \text{else},
            \end{cases}
    \end{array}
\end{equation}
where $e(\pi)$ is the Kottwitz sign for the pure inner form $\pi$ belongs to. This is Hypothesis \ref{hypothesis: KLH}.


Conversely, the equality \eqref{eq ind res intro} implies the above two expectations. Here is our main theorem.

\begin{thm}[{Propositions \ref{proposition: pure parabolic induction and restriction}, \ref{theorem: ind res adjunction implies klh}}]\label{theorem: main theorem in the introduction}
    Let $G$ be a quasisplit connected reductive group over $F$. Assume the enhanced local Langlands correspondence for $G$ and for its standard Levi subgroups, together with the following desiderata:
    \begin{enumerate}
        \item compatibility with the Langlands classification, Desideratum \ref{desideratum: Langlands classification};
        \item compatibility with isomorphisms, Desideratum \ref{desideratum: isomorphism};
        \item compatibility of parabolic induction with infinitesimal parameters, Desideratum \ref{desideratum: Haines}.
    \end{enumerate}
    Then the following two statements are equivalent.
    \begin{enumerate}
        \item[(a)] The induction of tempered representations, Desideratum \ref{desideratum: induction of tempered representations}, and the $p$-adic analogue of the Kazhdan--Lusztig hypothesis, Hypothesis \ref{hypothesis: KLH}, both hold for $G$ and for each of its standard Levi subgroups.
        \item[(b)] For every standard Levi subgroup $L \leq M \leq G$ with $\xi_L^M: \dualgroup{L}\hookrightarrow \dualgroup{M}$, every infinitesimal parameter $\lambda_L$ of $L$ with $\lambda_M=\xi_{L}^M\circ\lambda_L$, all $\pi\in K\Pi_{\lambda_L}^\pure(L)$ and all $\mathcal{F} \in K\Per_{H_{\lambda_M}}(V_{\lambda_M})$,
        \[\langle \Ind_{L}^M(\pi),\mathcal{F}\rangle = \langle \pi,(\xi_{L}^M)^*(\mathcal{F})\rangle.\]
        In other words, $\Ind_L^M=\prescript{t}{}{(\xi_{L}^M)^*}$ with respect to the pairing \eqref{equation: pairing in the introduction}. 
    \end{enumerate}
\end{thm}

\subsection{Outline of the proof}
It suffices to treat the case of a Levi subgroup $M$ of $G$ with $\xi: \Lgroup{}{M}\hookrightarrow \Lgroup{G}$. For (a)$\implies$(b) in this case, we wish to show for all $\mathcal{F}\in K\Per_{H_\lambda}(V_\lambda)$
\[\langle \Ind_M^G(\pi), \mathcal{F}\rangle = \langle \pi, \xi^*(\mathcal{F})\rangle.\] Both sides of the identity 
are linear, so it suffices to compare them on bases of standard objects in  $K\Pi_{\lambda_M}^\pure(M)$ and $K\Per_{H_\lambda}(V_\lambda)$. On right hand side, we use the basis of standard sheaves $\stsh{\tau}{C}$ of $K\Per_{H_\lambda}(V_\lambda)$, on which the restriction is computed in \cite{ABV92}*{Proposition 23.7(b)}; see \eqref{eqn: restriction formula on standards}:
\[
    \xi^*(\stsh{\tau}{C}) = \sum_{d\xi(D)\subseteq C,~\sigma \in \irrepdual{A}_D} m(\sigma,\tau) \cdot  \stsh{\sigma}{D}.
\]
Here, $D$ runs over the $H_{\lambda_M}$-orbits in $V_{\lambda_M}$ and $m(\sigma,\tau)$ is the multiplicity of $\sigma$ in $\tau \vert_{A_D}$. On the left hand side, we decompose $\Ind_M^G(\pi)$ using the basis of standard modules. From the decomposition \eqref{eq temp decomp intro}, we prove in Theorem \ref{theorem: induction formula} that parabolic induction decomposes standard modules in exactly the same way in the Grothendieck group:
\[\Ind_M^G(\Delta(\phi_M,\tau_M))=\sum_{\tau \in \irrepdual{A}_\phi} m(\tau_M,\tau)\, \Delta(\phi,\tau).\]
Then the $p$-adic analogue of the Kazhdan-Lusztig hypothesis says precisely that the pairing is diagonal on these two bases up to Kottwitz signs, which gives us the claim.


For the converse, let $\Ind_M^G(\pi_M)$ be either {a tempered induction or a standard module}. In both cases, the representation $\pi_M$ has an $L$-parameter $\phi_M$ whose corresponding $H_{\lambda_M}$-orbit in $V_{\lambda_M}$ is open, and the question becomes the evaluation of 
\[\langle \Ind_M^G(\pi_M), \stsh{\tau}{C_{\phi}} \rangle.\]
where $\phi=\xi \circ \phi_M$  and $\tau_{C_\phi}$ is the extension by zero of the local system corresponding to $\tau$ on $C_\phi$. Assuming \eqref{eq ind res intro}, this is equal to 
\[  \langle \pi_M,  \xi^{\ast} \stsh{\tau}{C_{\phi}} \rangle. \]
 Then again using \cite{ABV92}*{Proposition 23.7(b)}, the desired conclusion follows.

\subsection{Related work and applications}
Theorem \ref{theorem: main theorem in the introduction} is a compatibility statement. It is the Levi case of the ``Langlands functoriality map'' of \cite{ABV92}*{Definition 1.36}, which Adams, Barbasch and Vogan use to describe the expected behaviour of their geometric packets in the archimedean case. The present result is the $p$-adic counterpart for Levi subgroups, where it has a clean interpretation in terms of parabolic induction.

For $G=\GL_n$, the identity of Theorem \ref{theorem: main theorem in the introduction}(b) was proved in \cite{CunninghamRayVoganGLn}. In this case the component groups $A_\phi=1$, $H^1(F,\GL_n)=1$, and parabolic induction from a tempered representation of a Levi subgroup is irreducible, so that Desideratum \ref{desideratum: induction of tempered representations} holds automatically. The $p$-adic Kazhdan--Lusztig hypothesis for $\GL_n$ (which was long regarded as a folklore theorem) is proved by Solleveld \cite{Sol25} for substantial classes of groups and representations including all classical groups; see also the exposition of Balodis in \cite{balodis2026proofpadickazhdanlusztighypothesis}. The identity is a key input in the proof of Vogan's conjecture for $p$-adic $\GL_n$ in \cite{CunninghamRayVoganGLn}: for an Arthur parameter $\psi$ of $\GL_n$, the Arthur packet and the ABV-packet coincide, $\Pi_\psi(G)=\Pi^\ABV_{\phi_\psi}(G)$, where $\phi_\psi$ is the Langlands parameter attached to $\psi$. This conjecture appears in \cite{ABV92} for archimedean groups and in \cite{Vog93} for $p$-adic groups, and it is formulated in the language of Vogan varieties and vanishing cycles in \cite{CFMMX22}*{Conjecture 1}.

Theorem \ref{theorem: main theorem in the introduction}(b) is related to geometric realizations of the Langlands-Shelstad and Kottwitz-Shelstad transfers (that is, maps of Grothendieck groups of perverse sheaves whose transposes give these transfers). These geometric realizations constitute a key step in proving Vogan's conjecture for groups for which endoscopic classification is known. This is taken up in \cite{CHLLRX26}. Specifically, see Theorem 15.3 of \textit{loc. cit.}; an important difference being that there the `transfer' and `lift' maps originate from the Grothendieck group of stable distributions on a quasisplit endoscopic group, whereas in this paper we work with the full Grothendieck group and do not restrict ourselves to the quasisplit case of Levi subgroup. 

Combined with the Fourier transform on $V_\lambda$ as in \cite{CFMMX22}*{Section 10.2.4}, our result is a step towards a geometric description of the Aubert-Zelevinsky involution for general connected reductive $p$-adic groups: the involution is compatible with parabolic induction, and the Fourier transform is compatible with the restriction $\xi^*$. The identity can also be used to obtain a neat characterisation of perverse sheaves corresponding to supercuspidal representations. Indeed, let $\mathcal{P}(\pi) \in \Per_{H_\lambda}(V_\lambda)_{/\text{iso}}^{\mathrm{simple}}$. Then $\pi$ is supercuspidal if and only if $\xi^*\mathcal{P}(\pi)=0$ for all  proper Levi subgroups $M <G$ for which $\lambda$ factors through the embedding $\xi:{}^L M \hookrightarrow {}^L G$ as $\lambda_M$, and $\xi^*:K\Per_{H_\lambda}(V_\lambda)\to K\Per_{H_{\lambda_M}}(V_{\lambda_M})$ is the corresponding equivariant pullback. This is treated in Section \ref{section: supercuspidal sheaves}.

\subsection{Acknowledgments} 
The authors would like to thank Clifton Cunningham for his constant support.
The second named author would like to thank Julia Gordon at University of British Columbia for providing a conducive research environment during which this paper was written.
\subsection{AI declaration}
We used LLMs to edit the document for typos, and as a search and summarization tool. All the ideation and mathematical work in this document is done by humans.

\section{Background from representation theory}
Throughout, $F$ is a nonarchimedean local field of characteristic $0$ with residue field of size $q=p^f$, for a fixed prime $p$ and an integer $f$.
\subsection{Groups and representations}\label{Groups and representations}
Let $G$ be a quasisplit connected reductive algebraic group. Fix a Borel subgroup $B \subset G$ and a maximal torus $T \subseteq B$. 
Let \[X^*(T)=\operatorname{Hom}(T,\mathbb{G}_m),\qquad X_*(T)=\operatorname{Hom}(\mathbb{G}_m,T)\] denote the character and cocharacter lattices of $T$ respectively. 
Let $\Sigma \subset X^*(T)$
be the root system for $(G,T)$ and $\dualgroup{\Sigma} \subset X_*(T)$
the corresponding coroot system. The choice of $B$ determines a set $\Sigma^+$ of positive roots
and a corresponding set $\Delta$ of simple roots. For each simple root $\alpha\in\Delta$, let
$
X_\alpha:\mathbb{G}_a\longrightarrow U_\alpha
$
be a fixed isomorphism onto the corresponding root subgroup. The data
$
(G,B,T,\{X_\alpha\}_{\alpha\in\Delta})
$
is called a \emph{pinning} of $G$. The pinning determines a canonical splitting of the based root datum
$
(X^*(T),\Sigma,\Delta,X_*(T),\dualgroup{\Sigma},\dualgroup{\Delta}),
$
and identifies the group of pinned automorphisms of $G$ with the automorphism group of the based root datum. Throughout this paper, all constructions are understood with respect to this fixed pinning. Moreover, since we fix a Borel subgroup of $G$, we also get a collection of standard parabolic subgroups of $G$. Associated to the pinned group $(G,B,T,\{X_\alpha\})$ is the pinned Langlands dual group
$(\dualgroup{G},\dualgroup{B},\dualgroup{T},\{X^\vee_{\alpha^\vee}\}),
$ whose based root datum is dual to that of $G$. In this paper, we consider the complex points $\dualgroup{G}(\CC)$ of the dual group of $G$. {We often shorten the notation of the complex points of the dual group to simply $\dualgroup{G}$.}

Next we fix notation for the Weil group. Fix an algebraic closure $\overline{F}$ of $F$. The Weil group is a dense subgroup of the absolute Galois group $\Gal(\overline{F}/F)$ that satisfies the short exact sequence
\[1 \rightarrow I_F\rightarrow W_F \rightarrow \mathbb{Z} \rightarrow 1\]
where $I_F$ is the inertia subgroup and
$\Z$ is generated by the geometric Frobenius map $\Fr$ of $F$. We let $|\cdot|$ denote the absolute value on $W_F$ such that $|\Fr|=q$.

 The action of the absolute Galois group
$
\Gamma_F=\operatorname{Gal}(\overline{F}/F)
$
on the based root datum of $G$ induces an action on $\dualgroup{G}$ by pinned automorphisms. The $L$-group of $G$ is the semidirect product
\[
{}^LG = \dualgroup{G}\rtimes W_F.
\]
For any standard Levi subgroup $M$ of $G$, there is an embedding between their $L$-groups
\begin{align}\label{eq Levi L embedding}
    \xi: {}^L M \hookrightarrow {}^L G,
\end{align} well-defined up to conjugation by $\dualgroup{G}$. See \cite{Bor79}*{Section 15.1} or \cite{Tai25}*{Section 4.3} for details.

Finally, we fix notation for the two sets of representations that occur throughout. For a finite group $A$, let $\irrepdual{A}$ denote the set of equivalence classes of complex irreducible representations of $A$. 
We let $\Pi(G)$ denote the set of equivalence classes of complex smooth irreducible representations of $G(F)$. We let $K\Pi(G)$ denote the Grothendieck group of smooth representations of finite length of $G$.

\subsection{Parabolic induction}\rev{\label{subsec: parabolic induction}}
In this section, $G$ is a connected reductive algebraic group, not necessarily quasisplit.  For a parabolic subgroup $P$ of $G$ with a Levi subgroup $M$, we let $\Ind_P^G$ denote the normalized parabolic induction from finite length representations of $M(F)$ to representations of $G(F)$. More precisely, let $(\sigma,W)$ be an admissible representation of $M$ which is inflated to a representation of $P$. Then $\Ind_P^G(\sigma)=(\pi,V)$, where $V$ is the $\CC$-vector space of smooth functions $f:G(F)\to W$ such that 
\[f(pg)=\delta_P(p)^{1/2}\sigma(p)f(g) ~\text{ for }~ p\in P(F) ~\text{ and }~g\in G(F),\]
where $\delta_P$ is the modular quasicharacter of $P$. The $G$-action on $V$ is given by $\pi(g)f(x)=f(xg)$ for all $g,x\in G(F)$. 

Suppose that $M$ is also a Levi factor of another parabolic subgroup $P'$ of $G$. It is known that $\Ind_P^G(\sigma)$ and $\Ind_{P'}^G(\sigma)$ have the same image in $K\Pi(G)$ (\cite{BZ80}*{Theorem 2.9, Remark 2.10}). Therefore, we let 
\[ \Ind_M^G:K\Pi(M) \to K\Pi(G) \]
denote the linear transformation  $\sigma \mapsto \Ind_P^G(\sigma)$ for any choice of $P$.
\subsection{The enhanced local Langlands correspondence}\label{subsection: the enhanced local Langlands correspondence}
A Langlands parameter for $G$ is a continuous homomorphism $\phi: W_F \times \SL_2(\BC) \to \Lgroup{G}$ satisfying conditions specified in \cite{CFMMX22}*{Section 3.4}. An $L$-parameter is a $\dualgroup{G}$-conjugacy class of a Langlands parameter. We write $\phi$ for both the Langlands parameter and the $L$-parameter, unless it is necessary to emphasize the distinction. We let $\Phi(G)$ denote the set of $L$-parameters of $G$. 

In its classical form, the local Langlands correspondence is a conjectural many-to-one surjection between the set $\Pi(G)$ of smooth irreducible representations of $G(F)$ and $L$-parameters of $G$:
 \[ \LLC: \Pi(G) \to \Phi(G).\]
The fibres of $\LLC$ are the $L$-packets. In our setting, the conjectural surjection is refined to a bijection using the component group for each $\phi$,
\begin{equation}\label{eqn: a phi}
    A_\phi = Z_{\dualgroup{G}}(\phi)/Z_{\dualgroup{G}}(\phi)^0,
\end{equation}
where $Z_{\dualgroup{G}}(\phi)$ denotes the centralizer of the image of $\phi$ in the dual group $\dualgroup{G}$, and $Z_{\dualgroup{G}}(\phi)^0$ is its identity component. 
Let $\Phi^\mathfrak{e}(G)$ denote the set of $\dualgroup{G}$-conjugacy classes of pairs $(\phi, \tau)$, where $\phi$ is a Langlands parameter for $G$ and $\tau \in \irrepdual{A}_\phi$. Such a pair is called an \emph{enhanced $L$-parameter}. If we restrict our attention to the quasisplit group $G$ alone, the fiber of $\LLC$ over $\phi$ is parametrized by the quotient $A_{\phi}/Z(\dualgroup{G})^{\Gamma_F}$. However, to incorporate Vogan's geometric perspective (see Section~\ref{sec: geometric LLC}), it is more natural to work with $A_{\phi}$ itself. This viewpoint naturally leads us to consider the pure inner forms of $G$ simultaneously, which are parametrized by the Galois cohomology group $H^1(F,G)$. Let
\[\Pi^\pure(G)=\bigsqcup_{\delta\in H^1(F,G)}\Pi(G_\delta)\]
denote the set of equivalence classes of smooth irreducible representations of the quasisplit group $G$, which corresponds to the trivial element in $ H^1(F,G)$,
along with its other pure inner forms $G_{\delta}$. The enhanced local Langlands correspondence is a canonical bijection 
\begin{equation}\label{eqn: enhanced llc}
    \begin{array}{cc}
             \LLC^\mathfrak{e}: \Pi^\pure(G) &\to \Phi^\mathfrak{e}(G),  \\
         \pi&\mapsto (\phi, \tau). 
    \end{array}
\end{equation}
We often write $\pi(\phi,\tau)$ to denote the representation $\pi$ which satisfies $\LLC^\mathfrak{e}(\pi)=(\phi,\tau)$. 
We let $\LLC$ denote the composition of $\LLC^\mathfrak{e}$ with the projection $(\phi, \tau) \mapsto \phi$. Then we define $L$-packets
\[ \Pi_{\phi}^{\pure} (G):= \LLC^{-1}(\phi),\ \ \Pi_{\phi} (G_{\delta}):= \Pi_{\phi}^{\pure} (G) \cap \Pi(G_{\delta}).\]
The bijection \eqref{eqn: enhanced llc} is expected to depend on a choice of Whittaker datum $\mathfrak{w}$ of $G$, which we have fixed. We write $\LLC^{\mfr{e}}_{\mfr{w}}$ for $\LLC^{\mfr{e}}$ and $\pi_{\mfr{w}}(\phi,\tau)$ for $\pi(\phi,\tau)$ to emphasize the dependence on $\mfr{w}$ if necessary. 

The class $\delta$ attached to an enhancement is read off by the Kottwitz isomorphism, which we now recall. For a Langlands parameter $\phi$, the inclusion
\[Z(\dualgroup{G})^{\Gamma_F} \hookrightarrow Z_{\dualgroup{G}}(\phi)\] { induces the well-defined map of their component groups}
\begin{equation}\label{equation: inclusion of componenet grp of center in a phi}
    \iota_\phi: \pi_0(Z(\dualgroup{G})^{\Gamma_F}) \to A_\phi.
\end{equation}
Moreover, the image of $\iota_\phi$ is contained in the center of $A_\phi$. For any $\tau \in \irrepdual{A}_\phi$, the pullback $\tau \circ \iota_\phi$ is an isotypic direct sum of the pullback of the central character of $\tau$ to $\pi_0(Z(\dualgroup{G})^{\Gamma_F})$. We let $\chi_{\tau} \in \Hom(\pi_0(Z(\dualgroup{G})^{\Gamma_F}), \BC^\times)$ be the corresponding character. Recall the Kottwitz isomorphism for $G$:
\begin{equation}\label{eqn: kottwitz map}
    \kappa_G: \Hom(\pi_0(Z(\dualgroup{G})^{\Gamma_F}), \BC^\times) \xrightarrow[]{\sim} H^1(F,G).
\end{equation}
For $\tau \in \irrepdual{A}_\phi$, define 
\begin{equation}\label{eqn: kappa tau}
    \kappa_G(\tau):=\kappa_G(\chi_{\tau}).
\end{equation}

Often, we drop the subscript and use $\kappa(\tau)$ instead of $\kappa_G(\tau)$ when the group is clear. Finally, we build the following requirement into our definition of $\LLC^{\mathfrak{e}}$: for $(\phi,\tau)\in \Phi^{\mathfrak{e}}(G)$ and $\delta \in H^1(F,G)$,
\begin{equation}\label{equation: kappa determines the form}
    \kappa(\tau)=\delta \quad \Longleftrightarrow \quad \pi(\phi,\tau)\in \Pi(G_\delta).
\end{equation} We recall the definition of the Kottwitz sign.
Write $r(H)$ for the $F$-rank of the derived group of $H$. For $\delta \in H^1(F,G)$, the \emph{Kottwitz sign} of the pure inner form $G_\delta$ is
\begin{equation}\label{equation: kottwitz sign}
    e(\delta):= e(G_\delta)=(-1)^{r(G)-r(G_\delta)}.
\end{equation}
In this way, to each $\delta \in H^1(F,G)$, we assign the Kottwitz sign $\pm 1$. This may be interpreted as a function $e:H^1(F,G)\to \{\pm 1\}$. We can attach a Kottwitz sign to $\tau \in \Irrep(A_\phi)$ using \eqref{eqn: kappa tau}:
\begin{equation}\label{equation: kottwitz sign of tau}
    e(\tau):=e(\kappa_G(\tau)).
\end{equation}

With the requirement \eqref{equation: kappa determines the form} the enhanced local Langlands correspondence \eqref{eqn: enhanced llc} restricts to a bijection for each $\delta \in H^1(F,G)$: 
\begin{equation}\label{equation: llc delta}
    \LLC^{\mathfrak{e}}_{\delta}: \Pi(G_{\delta}) \to \{ (\phi, \tau) \in \Phi^\mathfrak{e}(G) \ \vert\  \kappa(\tau)=\delta \}.
\end{equation}

\subsection{Langlands classification}\label{subsec: Langlands classification}
We recall the Langlands classification, for representations and for $L$-parameters, together with its refinement to enhanced $L$-parameters.

Let $\delta\in H^1(F,G)$. We let $\Phi(G_{\delta})$ denote the subset of $\Phi(G)$ that is $G_{\delta}$-relevant. Namely, $\phi \in \Phi(G_{\delta})$ if and only if the minimal Levi subgroup of ${}^LG$ containing the image of $\phi$ is conjugate to the $L$-group of a standard Levi subgroup of $G_{\delta}$.

The Langlands classification for both $\Pi(G_{\delta})$ (\cite{Kon03}*{Theorem 3.5}) and $\Phi(G_{\delta})$ (\cite{SZ18}) provides canonical bijections
\begin{align*}
    \pi \leftrightarrow (P, \pi_\temp, \nu), \qquad \phi \leftrightarrow (P, \phi_\temp, \nu),
\end{align*}
where:
\begin{itemize}
    \item $\pi \in \Pi(G_{\delta})$ and $\phi \in \Phi(G_{\delta})$, 
    \item $P$ is a standard parabolic subgroup of $G_{\delta}$ with standard Levi component $M$,
    \item $\nu$ is a positive real unramified character of $M$ lying in the positive Weyl chamber determined by $P$, and
    \item $\pi_\temp \in \Pi(M)$ and $\phi_\temp \in \Phi(M)$ are both tempered (i.e., $\phi_\temp|_{W_F}$ has a bounded image).
\end{itemize}


Under this correspondence, $\pi$ is realized as the unique irreducible quotient of the \emph{standard module} (or \emph{standard representation})
\[
\Delta(\pi):=\Ind_P^{G_{\delta}}(\pi_\temp \otimes \nu),
\]
and the triple $(P,\pi_\temp,\nu)$ is uniquely determined by $\pi$.

For $L$-parameters, the positive real unramified character $\nu$ determines an unramified cocycle $c_{\nu}:W_F \to Z(\dualgroup{M})$. Define the twist $\phi_\temp \otimes \nu$  by 
\[(\phi_\temp \otimes \nu)(w,x):= c_{\nu}(w) \phi_{\temp}(w,x); \]
see \cite{SZ18}*{Definition 4.5}.
The relation between $\phi$ and the triple $(P, \phi_\temp, \nu)$ is given by
\[
\phi = ({}^L M \hookrightarrow {}^L G) \circ (\phi_\temp \otimes \nu).
\]
Since $\phi$ is $G_{\delta}$-relevant, the Levi subgroup $M$ is a Levi subgroup of $G_{\delta}$.

We need the following ingredients from the proof of the Langlands classification for $\Phi(G_{\delta})$.
\begin{prop}[{\cite{SZ18}*{Lemma 5.2, Propositions 5.5, 6.1, equation (6) in \S 7.2}}]\label{prop: SZ18}
   Let $\phi$ be a Langlands parameter of $G$. Suppose that the corresponding $L$-parameter $[\phi]$ corresponds to $(P, [\phi_\temp], \nu)$ under the Langlands classification. Then the following holds.
   \begin{enumerate}
    \item The hyperbolic part $z_{\phi}$ of $\phi(\Fr)$ is independent of the choice of $\Fr$. Moreover, the character $\nu$ is determined by the $\dualgroup{G}$-conjugacy class of $z_{\phi}$.
    \item We have ${}^LM $ conjugate to $ Z_{{}^L G}(z_{\phi})$ under $\dualgroup{G}$.
    \item The homomorphism of component groups $A_{\phi_\temp}^{\phi}(\xi): A_{\phi_\temp} \to A_{\phi}$ induced from the $L$-embedding $\xi: \prescript{L}{}{M} \hookrightarrow \prescript{L}{}{G}$ is an isomorphism.
   \end{enumerate}
\end{prop}

By Part (3) in the above proposition, the Langlands classification for $\Phi(G)$ can be further refined to the Langlands classification for $\Phi^{\mathfrak{e}}(G)$. More precisely, we define
\[ (\phi, \tau) \leftrightarrow (P, (\phi_\temp, \tau_\temp), \nu),\]
if $\phi \leftrightarrow (P, \phi_\temp, \nu)$ and $\tau_\temp = \tau \circ A_{\phi_\temp}^{\phi}(\xi)$.

It is also useful to define the set
\begin{equation}\label{eqn: a phi delta}
    \irrepdual{A}_{\phi,\delta}:=\{\tau \in \irrepdual{A}_\phi\ | \  \kappa(\tau)=\delta\},
\end{equation}
which allows us to write
\[\irrepdual{A}_\phi=\bigsqcup_{\delta\in H^1(F,G)}\irrepdual{A}_{\phi, \delta}.\]

\section{Background from geometry}\label{sec: geometric LLC}
The main goal of this section is to outline the geometric objects that inform Vogan's geometric perspective on the local Langlands correspondence from \cite{Vog93} which of course depends on the existence of $\LLC^{\mathfrak{e}}$. We interpret this geometric perspective using a pairing on Grothendieck groups.

\subsection{A moduli space for Langlands parameters}\label{subsec: a moduli space for langlands parameters}
An \emph{infinitesimal parameter} for $G$ is a continuous homomorphism $\lambda: W_F \to \Lgroup{G}$ that is a section of $\Lgroup{G} \to W_F$ and whose image consists of semisimple elements in $\Lgroup{G}$; see \cite{CFMMX22}*{Section 4.1}. The map $\phi \mapsto \lambda_\phi$ is surjective, and finite-to-one modulo equivalence.  We now fix an infinitesimal parameter $\lambda$ of $G$. The Vogan variety for $\lambda$ is 
\begin{equation}\label{equation: vogan variety}
    V_\lambda=\{x \in \mathrm{Lie}(\dualgroup{G})~ \vert ~ \Ad(\lambda(w))x=|w|x, ~\forall w\in W_F \}.
\end{equation}
It carries the $\Ad$-action of the centralizer group
\begin{equation}\label{equation: h lambda}
    H_\lambda=Z_{\dualgroup{G}}(\lambda).
\end{equation}
The complex affine variety $V_\lambda$ is a prehomogeneous vector space stratified under the $H_\lambda$-action into finitely many orbits. 

\subsection{Simple Perverse Sheaves}\label{subsec: simple perverse sheaves}
For an $H_\lambda$-orbit $C$ in $V_\lambda$, choose a base point $x \in C$. The equivariant fundamental group of $C$ is the finite group $A_C=\pi_0(Z_{H_\lambda}(x))$. Let $\Per_{H_\lambda}(V_\lambda)$ denote the category of $H_\lambda$-equivariant perverse sheaves on $V_\lambda$. It is the heart of the perverse $t$-structure on the $H_\lambda$-equivariant bounded derived category $D_{H_\lambda}(V_\lambda)$ of constructible sheaves on $V_\lambda$. Equivalence classes of simple equivariant local systems on $C$ are parametrized by the set $\irrepdual{A}_C$ of equivalence classes of irreducible representations of the finite group $A_C$. Let $\locsys{\tau}{C}$ denote the local system on $C$ corresponding to $\tau \in \irrepdual{A}_C$. The set $\Per_{H_\lambda}(V_\lambda)_{\text{/iso}}^{\text{simple}}$ of simple objects up to equivalence in $\Per_{H_\lambda}(V_\lambda)$ is given by
\begin{equation}\label{eq: simple perverse sheaves classified}
    \Per_{H_\lambda}(V_\lambda)_{\text{/iso}}^{\text{simple}}=\{\mathcal{IC}(\locsys{\tau}{C})~\vert ~ H_\lambda \text{-orbit }C \subseteq V_\lambda, \tau \in \irrepdual{A}_C \}.
\end{equation}
Here $\mathcal{IC}(\locsys{\tau}{C})$ denotes 
the corresponding intersection cohomology complex: given the inclusions $C\xhookrightarrow{j}\ol{C}\xrightarrow[]{i}V_\lambda$, 
\begin{equation}\label{equation: ic definition}
    \mathcal{IC}(\locsys{\tau}{C}):=i_*j_{!*}\locsys{\tau}{C}[\dim C]
\end{equation}
where $j_{!*}$ is the intermediate extension functor. See \cite{CFMMX22}*{Chapter 4}, or for more general exposition on perverse sheaves and derived categories for use in representation theory, see \cite{Ach21}*{Sections  3.1, 3.3, 6.4}. 

\subsection{Standard sheaves}\label{subsec: standards sheaves} The $H_\lambda$-equivariant bounded derived category $D_{H_\lambda}(V_\lambda)$ has another heart, the category $\operatorname{Sh}_{H_\lambda}(V_\lambda)$ of $H_\lambda$-equivariant constructible sheaves on $V_\lambda$; see \cite{Ach21}*{Section 6.4, Exercise 6.4.3}. We explain the objects of this category. Let $C$ be an $H_\lambda$-orbit in $V_\lambda$. For every local system $\locsys{\tau}{C}$ supported on $C$, we obtain a constructible sheaf on $V_\lambda$ via the extension by zero functor (see \cite{Ach21}*{Section 1.3}):
\begin{equation}
    \stsh{\tau}{C} :=k_!\locsys{\tau}{C} \text{ where } k:C \hookrightarrow V_\lambda.
\end{equation} 
We refer to the sheaves $\stsh{\tau}{C}$ as standard sheaves. The simple objects up to equivalence are
\begin{equation}\label{eq: standard sheaves classified}
    \operatorname{Sh}_{H_\lambda}(V_\lambda)_{\text{/iso}}^{\text{simple}}=\{\stsh{\tau}{C}~\vert ~ H_\lambda \text{-orbit }C \subseteq V_\lambda, \tau \in \irrepdual{A}_C \}.
\end{equation}

\subsection{Vogan--Langlands correspondence}\label{Vogan--Langlands correspondence}
We recall Vogan's perspective on the local Langlands correspondence for the quasisplit connected reductive algebraic group $G$ and its pure inner forms, established in \cite{Vog93}. As explained in \cite{CFMMX22}*{Proposition 4.2}, every point on the variety $V_\lambda$ corresponds to a Langlands parameter with infinitesimal parameter $\lambda$. They are related by 
\[ \phi \mapsto x_{\phi}=d(\phi|_{\SL_2})\begin{pmatrix}
    0 & 1\\0& 0
\end{pmatrix} \in V_{\lambda}. \]
Every $H_\lambda$-orbit in $V_\lambda$ corresponds to a conjugacy class of such a Langlands parameter: there is a bijection $\phi \mapsto C_{\phi}$ between $L$-parameters of $G$ that share the infinitesimal parameter $\lambda$, and the $H_{\lambda}$-orbits $C_{\phi}$ in $V_{\lambda}$. 
From
\cite{CFMMX22}*{Lemma 4.4}, for each Langlands parameter $\phi$,
we have that \begin{equation}\label{eqn: equivalence of component group and eq fund group}
    A_{C_{\phi}}\simeq A_{\phi}.
\end{equation}
Recall that the definition of $A_{C_{\phi}}$ depends on a basepoint in $C_{\phi}$, which we choose to be $x_{\phi}$ for the above isomorphism.

Let $\Pi_{\lambda}^\pure(G)$ denote the subset of representations in $\Pi^\pure(G)$ whose Langlands parameters share the infinitesimal parameter $\lambda$. Let $\pi(\phi,\tau)$ be the representation satisfying $\LLC^{\mathfrak{e}}(\pi(\phi,\tau))=(\phi, \tau)$, where $\tau \in \irrepdual{A}_\phi$; see \eqref{eqn: enhanced llc}. Vogan's geometric perspective on the local Langlands correspondence can be summarized as the bijection
\begin{equation}\label{eqn: vogan langlands correspondence}
    \begin{array}{rcl}
    \Pi_{\lambda}^\pure(G) 
    &\leftrightarrow& \Per_{H_\lambda}(V_\lambda)_{\text{/iso}}^{\text{simple}} \\
    \pi(\phi, \tau) 
    &\mapsto&
    \mathcal{IC}(\locsys{\tau}{C_{\phi}})=:\mathcal{P}(\pi).
\end{array}
\end{equation}
\subsection{Grothendieck groups}\label{subsec: grothendieck groups}
Let $KD_{H_\lambda}(V_\lambda), K\Per_{H_\lambda}(V_\lambda), $ and $ K\operatorname{Sh}_{H_\lambda}(V_\lambda)$ denote the Grothendieck groups of the categories $D_{H_\lambda}(V_\lambda)$, $\Per_{H_\lambda}(V_\lambda)$, and $\operatorname{Sh}_{H_\lambda}(V_\lambda)$, respectively. One can interpret $K\Per_{H_\lambda}(V_\lambda)$ (resp.\ $K\operatorname{Sh}_{H_\lambda}(V_\lambda)$) as the $\Z$-span of the elements of the set \eqref{eq: simple perverse sheaves classified} (resp. \eqref{eq: standard sheaves classified}).

As in \cite{ABV92}*{Equation 7.10(a)}, we have the identification
\[K\Per_{H_\lambda}(V_\lambda)=K\operatorname{Sh}_{H_\lambda}(V_\lambda)=KD_{H_\lambda}(V_\lambda).\]
In particular, $K\Per_{H_\lambda}(V_\lambda)$ can be interpreted as a $\mathbb{Z}$-module (or, after a change of scalars, a $\BC$-module) with two bases: one consisting of simple perverse sheaves as in \eqref{eq: simple perverse sheaves classified} and the other consisting of standard constructible sheaves as in \eqref{eq: standard sheaves classified},
\begin{equation}\label{equation: gr grp perverse sheaves}
\begin{aligned}
    K\Per_{H_\lambda}(V_\lambda)&=\mathrm{Span}_{\mathbb{Z}}\{\mathcal{IC}(\locsys{\tau}{C}): C \text{ an }H_\lambda\text{-orbit in }V_\lambda, \tau \in \irrepdual{A}_C \}\\
    &=\mathrm{Span}_{\mathbb{Z}}\{\stsh{\tau}{C}: C \text{ an }H_\lambda\text{-orbit in }V_\lambda, \tau \in \irrepdual{A}_C \} {.}
\end{aligned}
\end{equation}
The orbits can be labeled using Langlands parameters by invoking the Vogan-Langlands correspondence. For example, recalling $\irrepdual{A}_{\phi,\delta}$ from \eqref{eqn: a phi delta}, we define
\begin{equation}\label{equation: kper subspace delta}
    K\Per_{H_\lambda}(V_\lambda)^\delta := \operatorname{Span}_{\mathbb{Z}} \{\mathcal{IC}(\locsys{\tau}{C_\phi}): \lambda_\phi=\lambda \text{ and } \tau \in  \irrepdual{A}_{\phi,\delta}\} {.}
\end{equation}
{We will see the representation theory counterpart of the two bases in \eqref{equation: gr grp perverse sheaves} in Section \ref{subsec: Kazhdan--Lusztig hypothesis}.}

We record a case of special interest that relates the extension by zero of a local system with its intersection cohomology complex. Suppose $C$ is an open $H_\lambda$-orbit in $V_\lambda$, that is, $\ol{C}=V_\lambda$. Keeping \eqref{equation: ic definition} in mind, we write $j: C \hookrightarrow{} V_\lambda$ and $\mathcal{IC}(\locsys{\tau}{C})=j_{!*}\locsys{\tau}{C}[\dim C]$. Note $\locsys{\tau}{C}[\dim C] \in \Per(C)$. We compare $\stsh{\tau}{C}=j_!\locsys{\tau}{C}$ and $\mathcal{IC}(\locsys{\tau}{C})$ in $K\Per_{H_\lambda}(V_\lambda)=KD_{H_\lambda}(V_\lambda)$: 
\[
j_!\locsys{\tau}{C}[\dim C]
=
\mathcal{IC}(\locsys{\tau}{C})
+
\sum_{\substack{D\subsetneq \ol{C}\\ \gamma \in \irrepdual{A}_D}}
m(D,\gamma)
\mathcal{IC}(\locsys{\gamma}{D}),
\]
where $m(D,\gamma)\in\mathbb{Z}$. One can deduce this relation using \cite{BBD82}*{p. 54-55}. Shifting by an integer $d$ is the same as multiplication by $(-1)^d$ in the Grothendieck group. We thus obtain that the standard sheaf $\stsh{\tau}{C}$, up to the important sign $(-1)^{\dim C}$, is a sum of the perverse sheaf $\mathcal{IC}(\locsys{\tau}{C})$ and sheaves supported on the boundary $V_\lambda\setminus C$:
\begin{equation}\label{equation: standard sheaf as a sum of perverse sheaf and boundary terms}
    (-1)^{\dim C} \stsh{\tau}{C}=\mathcal{IC}(\locsys{\tau}{C}) + \sum_{\substack{D\subsetneq \ol{C}\\ \gamma \in \irrepdual{A}_D}}
m(D,\gamma)
\mathcal{IC}(\locsys{\gamma}{D})
\end{equation}
for some integers $m(D,\gamma)$.
\subsection{Pairing of Grothendieck groups on the basis of simple objects}\label{subsec: Pairing simple}
For every infinitesimal parameter $\lambda :W_F \to \Lgroup{G}$, the Vogan-Langlands correspondence, as in \eqref{eqn: vogan langlands correspondence}, fixes a perfect pairing between Grothendieck groups
\[
    K\Pi^\pure_\lambda(G) \times K\Per_{H_\lambda}(V_\lambda) \to \mathbb{Z}
\]
defined by
\begin{equation}\label{eqn:pairing}
\langle \pi,\mathcal{P}\rangle =
\begin{cases}
    (-1)^{d(\pi)}e(\pi) & \mathcal{P} = \mathcal{P}(\pi),\\
    0 & \text{otherwise},
\end{cases}
\end{equation}
where $\LLC^{\mathfrak{e}}(\pi)=(\phi,\tau)$, $d(\pi) \coloneqq \dim C_{\phi}$, and $e(\pi)$ is the Kottwitz sign attached to $\pi$, that is $e(\pi):=e(\delta)$ (see \eqref{equation: kottwitz sign}) where $\pi \in \Pi(G_\delta)$. This pairing is equivalent to the one in \cite{Vog93}*{(8.11)'(a)}; subsequently it appears in \cite{CFMMX22}*{p. 77} using the language of Vogan varieties. 
\section{Desiderata}
We gather the list of desiderata of the enhanced local Langlands correspondence needed in this work.

\subsection{Compatibility with Langlands classification}\label{subsec: Compatibility with Langlands classification}
Recall from Section \ref{subsec: Langlands classification} the Langlands classification for representations and $L$-parameters. We consider not only the local Langlands correspondence for pure inner forms of $G$, but also the local Langlands correspondence for pure inner forms of all standard Levi subgroups of $G$. These correspondences are compatible via Langlands classification in the following sense.
\begin{desideratum}\label{desideratum: Langlands classification}
Suppose that $\LLC^{\mathfrak{e}}(\pi) = (\phi, \tau)$. Let 
\[\pi \leftrightarrow (P_{\pi}, \pi_{\mathrm{temp}}, \nu_{\pi}),\ \ \ (\phi, \tau) \leftrightarrow (P_{\phi}, (\phi_{\mathrm{temp}}, \tau_{\mathrm{temp}}), \nu_{\phi})\] 
under the Langlands classification. Then we have $P_{\pi} = P_{\phi}$ and $\nu_{\pi} = \nu_{\phi}$, and furthermore,
\[
    \LLC^{\mathfrak{e}}(\pi_{\mathrm{temp}}) = (\phi_{\mathrm{temp}}, \tau_{\mathrm{temp}}).
\]
\end{desideratum}

\begin{remark}\label{rmk LC tempered}
Desideratum \ref{desideratum: Langlands classification} implies that if $\pi$ is essentially tempered, then its $L$-parameter $\phi= \phi_{\mathrm{temp}} \otimes \nu_{\phi}$ is a twist of a tempered $L$-parameter $\phi_{ \mathrm{temp}}$. In particular, since $\phi_{\mathrm{temp}}$ corresponds to the open orbit in the associated Vogan variety by \cite{cunningham2024whittakernormalizationpadicabvpackets}*{Proposition 3.2}, so does $\phi$.
\end{remark}

\subsection{Induction of tempered representations}\label{subsec: Induction of tempered representations}
The second desideratum involves the decomposition of the parabolic induction $\Ind_P^G (\pi(\phi_M,\tau_M))$ when $\pi(\phi_M,\tau_M)$ is a tempered representation of the Levi subgroup $M$ of $G$. It is known that this induction is semisimple and consists of tempered representations of $G$. We describe the subquotients of this induction in terms of the enhanced $L$-parameters of $G$.

First, we introduce some notation. Note that we do not assume $\phi$ is tempered in this paragraph.
Let $\xi: \prescript{L}{}{M} \hookrightarrow \prescript{L}{}{G}$ be the $L$-embedding. We let $\phi := \xi \circ \phi_M$. The embedding $\xi$ induces a homomorphism between the component groups of $\phi_M$ and $\phi$:
\[ A_{\phi_M}^{\phi}(\xi): A_{\phi_M} \to A_{\phi}.\]
For any $\tau \in \irrepdual{A}_{\phi}$, we decompose the pullback of $\tau$ along $A_{\phi_M}^{\phi}(\xi)$ into a direct sum of irreducible representations of $A_{\phi_M}$ as
\begin{align}\label{equation: pullback of tau}
\tau \circ A_{\phi_M}^{\phi}(\xi) = \bigoplus_{\tau_M \in \irrepdual{A}_{\phi_M}} m(\tau_M, \tau)\tau_M, 
\end{align}
where $m(\tau_M, \tau) \in \Z_{\geq 0}$ represents the multiplicity of $\tau_M$ in this decomposition. We now use this to describe the decomposition of the tempered induction. 

\begin{desideratum}\label{desideratum: induction of tempered representations}
    Let $\phi_M$ be a tempered $L$-parameter of $M$ and $\tau_M \in \irrepdual{A}_{\phi_M}$. Then the parabolic induction $\Ind_P^G (\pi(\phi_M,\tau_M))$ decomposes as
    \begin{align}\label{eq decomp of tempered induction}
        \Ind_{P}^{G} (\pi(\phi_M,\tau_M))=\bigoplus_{\tau \in \irrepdual{A}_{\phi}} m(\tau_M, \tau) \pi(\phi, \tau).
    \end{align}
\end{desideratum}

\begin{remark} {\label{remark: tempered induction from ECI}}
    A stronger desideratum is to require that the tempered $L$-packet satisfies endoscopic character identities $($see \cite{Kal16}*{Conjecture E}$)$. Indeed, by the compatibility of endoscopic transfer and parabolic induction $($\cite{Hir04}*{Lemma 6.1}$)$, the desired equation \eqref{eq decomp of tempered induction} can be obtained by a straightforward computation using the Fourier inversion formula. A detailed argument will be provided in \cite{KLT}.
\end{remark}

In Theorem \ref{theorem: induction formula} we extend the identity \eqref{eq decomp of tempered induction} to the standard modules of nontempered representations in the Grothendieck group.

\subsection{Compatibility with isomorphism}\label{subsection: compatibility with isomorphism}
 We fix an isomorphism of quasisplit algebraic groups $f: G_1 \to G_2$ defined over $F$ throughout this subsection. It induces a bijection $ H^1(F,G_1) \to H^1(F,G_2)$ and we identify these two sets. In particular, for any $\delta \in H^1(F,G_1)$, the isomorphism $f$ also induces an isomorphism
$f_{\delta}: G_{1,\delta} \simeq G_{2, \delta}$. Moreover, we also have a bijection
\begin{align*}
    \Pi(G_{2,\delta}) & \to \Pi(G_{1,\delta}),\\
    \pi_2 & \mapsto \pi_2 \circ f_{\delta}.
\end{align*}

 Let ${}^L f: {}^L G_1 \to {}^L G_2$ be the corresponding $L$-isomorphism. It induces a bijection $\Phi(G_1) \to \Phi(G_2)$ given by $\phi_1 \mapsto {}^L f \circ \phi_1$, and also an isomorphism $ A_{\phi_1} \simeq A_{{}^L f \circ \phi_1}$, which we also denote by ${}^L f$. Thus, we obtain a map between enhanced $L$-parameters
 \begin{align*}
     \Phi^{\mfr{e}}(G_{1}) &\to  \Phi^{\mfr{e}}(G_2),\\
     (\phi_1,\tau_1) & \mapsto ({}^L f \circ \phi_1, \tau_1 \circ ({}^Lf)^{-1} )=: {}^Lf (\phi_1,\tau_1).
 \end{align*}

We need a desideratum that the enhanced local Langlands correspondence on $G_1$ and $G_2$ is compatible with the above bijections. 
\begin{desideratum}\label{desideratum: isomorphism}
    
   Fix a Whittaker datum $\mfr{w}_2=(B_2, \chi_2)$ of $G_2$ and let $\mfr{w}_1:= (f^{-1}(B_2), \chi_2 \circ f)$ be the corresponding Whittaker datum of $G_1$. The enhanced local Langlands correspondences for $G_1$ and $G_2$ satisfy the following: for all $(\phi_1, \tau_1) \in \Phi^{\mfr{e}}(G_1)$ with $\kappa_{G_1}(\tau_1)= \delta$, we have 
    \[ \pi_{\mfr{w}_1}(\phi_1,\tau_1)= \pi_{\mfr{w}_2}({}^Lf (\phi_1,\tau_1))   \circ f_{\delta}. \]
\end{desideratum}

We now describe the setting in which this desideratum will be applied in the
later section. Write $A_0$ for the maximal split subtorus of the maximal torus
$T$ of $G$ fixed in Section \ref{Groups and representations}. Let $M$ be a
standard Levi subgroup of $G$ and let $g \in N_G(A_0)(F)$ be an element such
that $M' := g M g^{-1}$ is again a standard Levi subgroup. Conjugation by $g$
then defines an isomorphism $f \colon M \to M'$, and the induced
$L$-isomorphism ${}^L f \colon {}^L M \to {}^L M'$ is given by conjugation by
some $\widehat{g} \in G^{\vee}$. For $\delta \in H^1(F,M)$ and
$\pi \in \Pi(M_{\delta})$ we write $\pi^g := \pi \circ f_{\delta}^{-1}$, where
$f_{\delta} \colon M_{\delta} \to M'_{\delta}$ is the isomorphism of pure inner
forms induced by $f$. For $(\phi, \tau) \in \Phi^{\mfr{e}}(M)$ we write
$({}^{\widehat{g}}\phi, {}^{\widehat{g}}\tau) := {}^L f (\phi,\tau)$.

Fix a Whittaker datum $\mfr{w} = (B,\chi)$ for $G$, and let $\mfr{w}_{M}$ and
$\mfr{w}_{M'}$ be the Whittaker data for $M$ and $M'$ obtained by restriction. Namely, write $N_B$ the unipotent radical of the Borel subgroup $B$, we have
\[ \mfr{w}_M = \bigl( B \cap M , \chi|_{N_B \cap M} \bigr), \]
and similarly for $M'$. Note that as long as $B$ is contained in a parabolic subgroup $P$ with Levi component $M$, the $M(F)$-conjugacy class of this restriction only depends on the $G(F)$-conjugacy class of $\mfr{w}$.

Now we check that $\mfr{w}_M$ and $\mfr{w}_{M'}$ are
related as in Desideratum \ref{desideratum: isomorphism}. Namely, the
pullback
\[ f^{*}\mfr{w}_{M'} := \bigl( f^{-1}(B \cap M'),\;
   \chi|_{N_B \cap M'} \circ f \bigr) \]
is equivalent to $\mfr{w}_M$. Indeed, $f^{*}\mfr{w}_{M'}$ is the restriction to
$M$ of the Whittaker datum $\mfr{w}' := (g^{-1}Bg,\; \chi \circ \Ad(g))$ of
$G$. Since $\mfr{w}'$ and  $\mfr{w}$ are conjugate by
$g \in G(F)$, their restriction is conjugate by $M(F).$
Hence, $\mfr{w}'_M = f^{*}\mfr{w}_{M'}$ is equivalent to
$\mfr{w}_M$. The desideratum therefore yields
\[ \LLC^{\mfr{e}}_{\mfr{w}_{M'}}
   \bigl( \pi_{\mfr{w}_M}(\phi,\tau)^g \bigr)
   = ({}^{\widehat{g}} \phi, {}^{\widehat{g}} \tau) \]
for every $(\phi, \tau) \in \Phi^{\mfr{e}}(M)$.

Finally, suppose we start instead from $\delta \in H^1(F,G)$, a standard Levi
subgroup $M_{\delta}$ of $G_{\delta}$, and an element $g_{\delta} \in
G_{\delta}(F)$ such that $M_{\delta}' := g_{\delta} M_{\delta} g_{\delta}^{-1}$
is again standard. Then one can choose $g \in N_G(A_0)(F)$ conjugating the
quasisplit Levi subgroup $M$ to $M'$ in such a way that the induced map
$f_{\delta} \colon M_{\delta} \to M_{\delta}'$ agrees with conjugation by
$g_{\delta}$. Then the discussion above applies.

\begin{remark}
    The description of $L$-parameter in Desideratum \ref{desideratum: isomorphism} is a special case of \cite{Bor79}*{Condition 10.3(5)}. The description of the enhanced $L$-parameters can be deduced from \cite{Kal23}*{Section 2.3.4}.
\end{remark}

\subsection{\texorpdfstring{Haines conjecture and the $p$-adic Kazhdan--Lusztig hypothesis}{}}\label{subsec: Kazhdan--Lusztig hypothesis}
In this subsection, we discuss the decomposition of the standard modules in terms of local Langlands correspondence.

First, we state Haines' conjecture on compatibility of induction and infinitesimal parameter \cite{Hai14}*{Conjecture 5.2.2}.
\begin{desideratum}\label{desideratum: Haines}
Let $\lambda_M$ be an infinitesimal parameter of $M$, a Levi subgroup of $G$. Let $\lambda:= ({}^L M \hookrightarrow {}^L G) \circ \lambda_M$ be the corresponding infinitesimal parameter of $G$. Suppose that $\pi_M  \in \Pi^{\pure}_{\lambda_M}(M)$ and $\pi$ is an irreducible subquotient of $\Ind_M^G\pi_M$. Then $ \pi \in \Pi^{\pure}_{\lambda}(G)$.
\end{desideratum}

An important special case of the above desideratum is that if $\pi \in \Pi_{\lambda}^{\pure}(G)$, then $\Delta(\pi) \in K\Pi_{\lambda}^{\pure}(G)$. Since the change of basis matrix from standard modules to irreducible representations is unitriangular (see \cite{BW00}*{Lemma XI.2.13}), the set 
\[ \{\Delta(\pi)\ | \ \pi \in \Pi^{\pure}_{\lambda}(G) \}\]
forms a $\Z$-basis. Therefore, assuming Desideratum \ref{desideratum: Haines}, the Grothendieck group $K\Pi_\lambda^\pure(G)$ has two bases {:} one consisting of equivalence classes of irreducible representations $\pi \in \Pi_\lambda^\pure(G)$ and the other consisting of equivalence classes of their standard modules $\Delta(\pi)$:
\begin{equation}\label{equation: gr grp of reps}
    K\Pi_\lambda^\pure(G)=\mathrm{Span}_{\mathbb{Z}}~\Pi_\lambda^\pure(G)=\mathrm{Span}_{\mathbb{Z}}\{\Delta(\pi)\ | \  \pi \in \Pi_\lambda^\pure(G)\}.
\end{equation}
Similarly, for $\delta\in H^1(F,G)$,
\[K\Pi_\lambda(G_\delta)=\mathrm{Span}_{\mathbb{Z}}~\Pi_\lambda(G_\delta)=\mathrm{Span}_{\mathbb{Z}}\{\Delta(\pi)\ | \  \pi \in \Pi_\lambda(G_{\delta})\}.\]

The Kazhdan--Lusztig hypothesis relates the change of basis matrix from irreducible representations to standard modules in $K\Pi_\lambda^\pure(G)$ to the change of basis matrix from simple perverse sheaves to standard constructible sheaves in $K\Per_{H_\lambda}(V_\lambda)$. See Section \ref{subsec: grothendieck groups} for aforementioned bases of Grothendieck groups. The statement of the hypothesis agrees with \cite{Sol25}*{Theorem E} and with \cite{Vog93}*{Conjecture 8.11} after dropping signs in the latter. This sign correction 
is illustrated via a simple example in \cite{balodis2026proofpadickazhdanlusztighypothesis}*{Example 3.5}. We state the hypothesis below in terms of the evaluation of the pairing \eqref{eqn:pairing} on the bases of standard objects in $K\Pi_\lambda^\pure(G)$ and $K\Per_{H_\lambda}(V_\lambda)$. The relation between this formulation and the change of basis matrices is explained in \cite{CunninghamRayVoganGLn}*{Section 1.4}.

\begin{hypothesis}[Kazhdan--Lusztig Hypothesis for $G$]\label{hypothesis: KLH}
Desideratum \ref{desideratum: Haines} holds. Moreover, the pairing defined in \eqref{eqn:pairing} is equivalent to the pairing defined by
\begin{equation}\label{equation: pairing on standards}
    \begin{array}{rcl}
        K\Pi_\lambda^\pure(G) \times K\Per_{H_\lambda}(V_\lambda)
            &\to& \mathbb{Z}\\
        \langle \Delta(\pi), \stsh{\tau}{C_{\phi}} \rangle 
            &\mapsto& 
            \begin{cases}
                e(\pi) & \LLC^{\mathfrak{e}}_G(\pi) = (\phi,\tau), \\
                0 & \text{else},
            \end{cases}
    \end{array}
\end{equation}
where $\stsh{\tau}{C_\phi}$ denotes the standard sheaf corresponding to $\tau$ on the $H_\lambda$-orbit $C_\phi$, and $e(\pi)$ is the Kottwitz sign of $\pi$ as in \eqref{eqn:pairing}. 
\end{hypothesis}

\section{Parabolic induction and restriction as maps between Grothendieck groups}
We recall and prove identities for parabolic induction of representations and equivariant restriction of sheaves at the level of Grothendieck groups.

\subsection{Setting}\label{subsection: setting}
Let $G$ be a connected reductive quasisplit algebraic group over a nonarchimedean local field $F$. Let $M$ be a standard Levi subgroup and $P$ be the standard parabolic subgroup with Levi component $M$. Recall from \eqref{eq Levi L embedding} that we have an $L$-embedding 
\begin{equation}\label{equation: restriction of L-groups}
    \xi: \prescript{L}{}{M} \hookrightarrow \prescript{L}{}{G}.
\end{equation} 


Let $\lambda$ be an infinitesimal parameter for $G$ which factors through an infinitesimal parameter $\lambda_M$ of $M$:
\begin{equation}\label{equation: setting inf par factors}
\begin{tikzcd}
	{W_F} && {\Lgroup{G}} \\
	& {\Lgroup{M}}
	\arrow["\lambda", from=1-1, to=1-3]
	\arrow["{\lambda_M}"', from=1-1, to=2-2]
	\arrow["\xi"', hook, from=2-2, to=1-3]
\end{tikzcd}
\end{equation}
The $L$-embedding \eqref{equation: restriction of L-groups} restricts to an
embedding of the complex algebraic groups:
\begin{equation}\label{equation: embedding of complex dual groups}
    \xi: \dualgroup{M} \hookrightarrow \dualgroup{G}.
\end{equation}

Let $M_\gamma$ (resp. $G_\delta$) correspond to the pure inner form $\gamma \in H^1(F,M)$ (resp. $\delta \in H^1(F,G)$). The inclusion $\xi: \dualgroup{M} \hookrightarrow \dualgroup{G}$ induces a homomorphism $H_\xi: H^1(F,M) \xrightarrow{} H^1(F,G) $, which we realize using Kottwitz isomorphism \eqref{eqn: kottwitz map} below. The inclusion $\xi$ induces the inclusion $Z(\dualgroup{G})^{\Gamma_F} \hookrightarrow Z(\dualgroup{M})^{\Gamma_F}$ which gives the map
    \[\pi_0(Z(\dualgroup{G})^{\Gamma_F}) \rightarrow \pi_0(Z(\dualgroup{M})^{\Gamma_F}).\]
    This in turn induces the map 
    \begin{align}\label{eq: Z xi}
        Z_\xi: \Hom(\pi_0(Z(\dualgroup{M})^{\Gamma_F}), \mathbb{C}^\times) \rightarrow \Hom(\pi_0(Z(\dualgroup{G})^{\Gamma_F}), \mathbb{C}^\times).
    \end{align}
    Then we obtain the following diagram.
\begin{equation}\label{equation: definition of h xi}
\begin{tikzcd}
	{H^1(F,M)} & {H^1(F,G)} \\
	{\operatorname{Hom}(\pi_0(Z(\dualgroup{M})^{\Gamma_F}), \mathbb{C}^\times)} & {\operatorname{Hom}(\pi_0(Z(\dualgroup{G})^{\Gamma_F}), \mathbb{C}^\times)}
	\arrow["{H_\xi}", dashed, from=1-1, to=1-2]
	\arrow["{\kappa_M^{-1}}"', from=1-1, to=2-1]
	\arrow["{Z_\xi}"', from=2-1, to=2-2]
	\arrow["{\kappa_G}"', from=2-2, to=1-2]
\end{tikzcd}
\end{equation}
Finally, it follows from \cite{Kot83}*{Proposition, part (6)} that
\begin{align}\label{eqn Kottwitz sign Levi}
    e(\gamma)=e(H_\xi(\gamma)).
\end{align}

\subsection{Restriction of sheaves}\label{subsec: restriction of sheaves}

Recall from \eqref{equation: h lambda} that $H_\lambda=Z_{\dualgroup{G}}(\lambda)$, so  $\xi$ from \eqref{equation: embedding of complex dual groups} restricts to the inclusion
\begin{equation}\label{equation: inclusion h lambda}
    \xi: H_{\lambda_M} \hookrightarrow H_\lambda.
\end{equation}
On the other hand, the map $\xi$ also induces a map $\operatorname{Lie}(\dualgroup{M})\to \operatorname{Lie}(\dualgroup{G})$ which restricts to an inclusion of Vogan varieties
\begin{equation}\label{equation: inclusion of vogan varieties}
    d\xi: V_{\lambda_M} \hookrightarrow{} V_\lambda.
\end{equation}
This map $d\xi$ is equivariant with respect to the $H_{\lambda_M}$-action on $V_{\lambda_M}$ and $H_\lambda$-action on $V_\lambda$:
\[d\xi(h\cdot v)=\xi(h)\cdot d\xi(v) \hspace{0.3cm} \forall h\in H_{\lambda_M}, v\in V_{\lambda_M}.\]
For any $x \in D \subset V_{\lambda_M}$, where $D$ is an $H_{\lambda_M}$-orbit in $V_{\lambda_M}$, the $H_\lambda$-orbit $C$ of $d\xi(x)\in V_\lambda$ is called the \emph{saturation} of $D$ in $V_\lambda$. From \eqref{equation: inclusion h lambda} we obtain the embedding $Z_{H_{\lambda_M}}(x) \hookrightarrow Z_{H_\lambda}(d\xi(x))$ which induces the map
\begin{equation}\label{equation: map of fundamental groups}
    A_D^C(\xi): A_D \to A_C
\end{equation}
between the equivariant fundamental groups of the orbits $D \subset V_{\lambda_M}$ and $C\subset V_\lambda$. Note the dependence on the choice of basepoint $x\in D$ (which determines $d\xi(x)\in C$). 

For a finite dimensional representation $\tau: A_C \to \GL(V_\tau)$, the pullback $\tau\vert_{A_D}$ is the representation $\tau \circ A_D^C(\xi)$.
\begin{equation}\label{equation: pullback of a rep of fund group}
\begin{tikzcd}
	{A_D} & {A_C} \\
	& {\GL(V_\tau)}
	\arrow["{A_D^C(\xi)}", from=1-1, to=1-2]
	\arrow["{\tau\vert_{A_D}}"', from=1-1, to=2-2]
	\arrow["\tau", from=1-2, to=2-2]
\end{tikzcd}
\end{equation}
The map \eqref{equation: inclusion of vogan varieties} also induces the equivariant pullback 
\[\xi^*:D_{H_\lambda}(V_\lambda) \to D_{H_{\lambda_M}}(V_{\lambda_M}) \]
between equivariant derived categories; see \cite{BL94} for the construction of the equivariant derived category and of the pullback functor along an equivariant map. 

Recall from Section \ref{subsec: grothendieck groups} that the Grothendieck groups coincide: $K\Per_{H_\lambda}(V_\lambda)=KD_{H_\lambda}(V_\lambda)$. This allows us to define the linear map which we also denote as $\xi^*$:
\begin{equation}\label{equation: restriction linear map perverse sheaves}
    \xi^*: K\Per_{H_\lambda}(V_\lambda) \xrightarrow{} K\Per_{H_{\lambda_M}}(V_{\lambda_M}).
\end{equation} 
Recall from the discussion in Section \ref{subsec: grothendieck groups} that the collection of standard sheaves $\stsh{\tau}{C}$ forms a basis of $K\textrm{Sh}_{H_{\lambda}}(V_{\lambda})= K \Per_{H_{\lambda}}(V_{\lambda})$. Therefore, the linear map $\xi^{\ast}$ is determined by its value on these standard sheaves, which is described explicitly in \cite{ABV92}*{Proposition 23.7(b)}. We recall the formula in the following proposition.

\begin{prop}\label{prop: equivariant restriction}
   In the setting above, we have
    \begin{equation}\label{eqn: restriction formula on standards}
    \xi^*(\stsh{\tau}{C}) = \sum_{\xi(D)\subset C,~\sigma \in \irrepdual{A}_D} m(\sigma,\tau) \cdot \stsh{\sigma}{D}.
\end{equation}
Here $D$ runs over the $H_{\lambda_M}$-orbits in $V_{\lambda_M}$ whose image under $d\xi$ is contained in $C$, $\sigma$ runs over $\irrepdual{A}_D$, and $m(\sigma,\tau)$ denotes the multiplicity of $\sigma$ in $\tau\vert_{A_D}$, in the notation of \eqref{equation: pullback of a rep of fund group}.
\end{prop}

\begin{remark}
   If $D=C_{\phi_M}$ and $C=C_\phi$ for Langlands parameters $\phi_M$ of $M$ and $\phi=\xi\circ\phi_M$ of $G$, then under the identifications \eqref{eqn: equivalence of component group and eq fund group} the map $A_D^C(\xi)$ is the map $A_{\phi_M}^\phi(\xi)$ of \eqref{equation: pullback of tau}. As a result, the multiplicities $m(\sigma,\tau)$ above agree with those of Section \ref{subsec: Induction of tempered representations}. 
\end{remark} 
\subsection{Lift}\label{subsection: lift}
For $\delta\in H^1(F,G)$, define the obvious inclusion and projection
\begin{equation}\label{equation: inclusion of grothendieck groups of reps}
    \iota_\delta^G: K\Pi_{\lambda}(G_\delta) \to K\Pi_\lambda^\pure(G), \ \ \pr_\delta^G:K\Pi_\lambda^\pure(G) \to K\Pi_{\lambda}(G_\delta). 
\end{equation}
Note that for any $\mathcal{F}\in K\Per_{H_\lambda}(V_\lambda)^\delta$,
\begin{equation}\label{equation: transpose of projection}
\langle \pr^G_\delta(\pi), \mathcal{F}\rangle = \langle \pi, \mathcal{F}\rangle.
\end{equation}
We often drop the superscript $G$ and write $\iota_\delta$ and $\pr_\delta$ when the group is clear.
With respect to the pairing in \eqref{eqn:pairing}, we have the transpose of pullback: 

\[\prescript{t}{}{\xi^*}: K\Pi_{\lambda_M}^\pure(M) \xrightarrow{} K\Pi^\pure_\lambda(G),\] characterized by the requirement that, for all $\mathcal{F} \in K\Per_{H_{\lambda}}(V_{\lambda})$ and $\pi \in K\Pi_{\lambda_M}^\pure(M) $, 
\begin{equation}\label{equation: transpose of restriction}
    \langle \prescript{t}{}{\xi^*}(\pi), \mathcal{F} \rangle = \langle \pi, \xi^*(\mathcal{F})\rangle.
\end{equation}

Now define
\begin{equation}\label{equation: lift}
 \Lift_{M_\gamma}^{G_\delta}:=\operatorname{pr}^G_\delta \circ \prescript{t}{}{\xi}^* \circ\iota^M_\gamma.
\end{equation}
That is,
\[\begin{tikzcd}
	{K\Pi_{\lambda_{M}}(M_\gamma)} & {K\Pi_\lambda(G_\delta)} \\
	{K\Pi_{\lambda_{M}}^{\pure}(M)} & {K\Pi^\pure_\lambda(G)}
	\arrow["{\Lift_{M_\gamma}^{G_\delta}}", from=1-1, to=1-2]
	\arrow["\iota_\gamma"', hook', from=1-1, to=2-1]
	\arrow["{\prescript{t}{}{\xi^*}}"', from=2-1, to=2-2]
	\arrow["{\operatorname{pr_\delta}}"', from=2-2, to=1-2]
\end{tikzcd}.\]

\subsection{Parabolic induction of representations}\label{subsec: parabolic induction of reps linear map}
Assume we are in the setting \ref{subsection: setting}. Recall from Section  {\ref{subsec: grothendieck groups}} the Grothendieck group $K\Pi_{\lambda}^\pure(G)$ and the parabolic induction functor $\Ind_P^G$ from Section \ref{subsec: parabolic induction}. Assuming Desideratum \ref{desideratum: Haines}, the linear map $\Ind_M^G$ restricts to
\begin{align}\label{equation: induction definition}
    \Ind_M^G: K\Pi^\pure_{\lambda_M}(M) &\xrightarrow{} K\Pi^\pure_{\lambda}(G),\\
    \pi &\mapsto \Ind_M^G(\pi). {\notag}
\end{align}
Here if $\pi \in K\Pi(M_{\gamma})$, then $\Ind_M^G(\pi)$ means $\Ind_{M_{\gamma}}^{G_{H_{\xi}(\gamma)}}(\pi)$ by abuse of notation.

Suppose $\phi$ is a Langlands parameter of $G$ that factors through the Langlands parameter $\phi_{M}$ of $M$ via $\xi$.
\begin{equation} {\label{equation: setting L par factors}}
\begin{tikzcd}
	{W_F\times \SL_2(\BC)} && {\Lgroup{G}} \\
	& {\Lgroup{M}}
	\arrow["\phi", from=1-1, to=1-3]
	\arrow["{\phi_M}"', from=1-1, to=2-2]
	\arrow["\xi"', hook, from=2-2, to=1-3]
\end{tikzcd}
\end{equation}
Recall from Section \ref{subsec: Induction of tempered representations} the homomorphism of component groups $A_{\phi_M}^\phi(\xi):A_{\phi_M}\to A_\phi$ induced by $\xi$, the pullback $\tau\vert_{A_{\phi_M}}:=\tau\circ A_{\phi_M}^\phi(\xi)$ of $\tau \in \irrepdual{A}_\phi$, and the multiplicities $m(\tau_M,\tau)$ of \eqref{equation: pullback of tau}.
Recall also the notation $\kappa_G(\tau)$ for $\tau \in \irrepdual{A}_\phi$ from \eqref{eqn: kappa tau} defined using the Kottwitz isomorphism and the map $H_\xi: H^1(F,M)\to H^1(F,G)$ from \eqref{equation: definition of h xi}. 
\begin{lemma}\label{lemma: relating galois colomology class of irreps of component groups}
     Suppose the irreducible representation $\tau_M\in \irrepdual{A}_{\phi_M}$ occurs in the pullback $\tau\vert_{A_{\phi_M}}$ according to \eqref{equation: pullback of tau}. Then, $\kappa_G(\tau)=H_\xi(\kappa_M(\tau_M))$.
\end{lemma}
\begin{proof}
    For $\tau \in \irrepdual{A}_\phi$, let $\chi_\tau \in \Hom(\pi_0(Z(\dualgroup{G})^{\Gamma_F}), \CC^\times)$ be the character in the pullback of $\tau$ via the map in $\iota_\phi$ in \eqref{equation: inclusion of componenet grp of center in a phi}. By definition \eqref{eqn: kappa tau}, $\kappa_G(\tau)=\kappa_G(\chi_\tau)$. Similarly $\kappa_M(\tau_M)=\kappa_M(\chi_{\tau_M})$. The conclusion follows from the diagram \eqref{equation: definition of h xi}.  
\end{proof}

Here is the main result of this section.

\begin{thm}\label{theorem: induction formula}
Assume that there exists an $\LLC^{\mathfrak{e}}$ for Levi subgroups of $G$ satisfying Desiderata \ref{desideratum: Langlands classification}, \ref{desideratum: induction of tempered representations}, \ref{desideratum: isomorphism}, and \ref{desideratum: Haines}. Recall the notation set up in Section \ref{subsec: Induction of tempered representations}. Suppose now that $\phi$ is \emph{any} $L$-parameter that factors through the parameter $\phi_M$ of a standard Levi subgroup $M$ of $G$. We have in the Grothendieck group $K\Pi_\lambda^\pure(G)$
\begin{equation}\label{induction formula}
   \Ind_{M}^{G} (\Delta(\phi_M,\tau_M))=\sum_{\tau \in \widehat{A}_{\phi}} m(\tau_M, \tau) \Delta(\phi, \tau).
\end{equation}
\end{thm}
\begin{proof}
    Write $\phi_M \leftrightarrow (P_0, \phi_{M_0}, \nu_0)$ as an $L$-parameter of $M$. Thus $P_0$ is a standard parabolic subgroup of $M$ with standard Levi component $M_0$, the parameter $\phi_{M_0}$ is tempered, and $\nu_0 \in \mfr{a}^*_{M_0}$ is positive with respect to $P_0$.  By Proposition \ref{prop: SZ18}, we can assume
    ${}^LM_0 = Z_{{}^L M}(z_{\phi_M})$, where $z_{\phi_M}$ is the hyperbolic part of $\phi_M(\Fr) \in {}^L{M}$ for any Frobenius element $\Fr$ of $F$. Let $\xi_{M_0}^M: \prescript{L}{}{M_0} \hookrightarrow \prescript{L}{}{M}$ be the $L$-embedding, and let $\tau_{M_0}:= \tau_M \circ A_{\phi_{M_0}}^{\phi_M}(\xi_{M_0}^M)$, which is irreducible by Proposition \ref{prop: SZ18}(3). In particular, $\tau_{M_0} \in \irrepdual{A}_{\phi_{M_0}}$. Let $P$ be the standard parabolic subgroup of $G$ with $P\cap M=P_0$, so that $M_0$ is the Levi component of $P$. Inducing in stages, the left hand side of \eqref{induction formula} is
    \[ \Ind_{M}^G \big(\Ind_{M_0}^M (\pi(\phi_{M_0},\tau_{M_0}) \otimes \nu_0)\big) = \Ind_{M_0}^G (\pi(\phi_{M_0},\tau_{M_0}) \otimes \nu_0).\]
    We show that the right-hand side is a sum of standard modules (with multiplicities) in the Grothendieck group.

    Although $\nu_0$ is positive with respect to $P_0$, it is not necessarily positive with respect to $P$. 
    Therefore, we take an element $w_1 \in N_G(A_0)(F)$, where $A_0$ is the maximal split torus of the maximal torus $T$ fixed in Section \ref{Groups and representations}, such that $w_1(\nu_0)$ lies in the closure of the Weyl chamber determined by the Borel subgroup we fixed. Then the character $w_1(\nu_0)$ of $w_1M_0w_1^{-1}$ extends to a character of a standard Levi subgroup $M_1:= Z_{G}(w_1(\nu_0))$, which is positive with respect to the corresponding standard parabolic subgroup $P_1$. Since $w_1M_0w_1^{-1}$ contains the centralizer of $A_0$, it is a semi-standard Levi subgroup of $M_1$. We may take $w_2 \in N_{M_1}(A_0)$ such that $M_0':= w_2w_1M_0w_1^{-1}w_2^{-1}$ is a standard Levi subgroup of $M_1$ (and of $G$). Set $w:=w_2w_1$. For a representation $\sigma$ of a pure inner form $M_{0,\gamma_0}$ of $M_0$, we let $\sigma^w$ denote the corresponding representation of the pure inner form $M_{0,w_{\ast}\gamma_0}'$ of $M_0'$. Note that the image of $\gamma_0 \in H^1(F,M_0) $ and $w_{\ast}\gamma_0 \in H^1(F,M_0')$ in  $ H^1(F,G)$ agree since $w \in G(F)$.
    Now we can continue to rewrite the left hand side of \eqref{induction formula} in Grothendieck group as  
\begin{align*}
     \Ind_{M_0}^G (\pi(\phi_{M_0},\tau_{M_0}) \otimes \nu_0)&=\Ind_{w M_0 w^{-1}}^G (\pi(\phi_{M_0},\tau_{M_0})^w \otimes w(\nu_0))\\
     &=\Ind_{M_1}^G (\Ind_{M_0'}^{M_1} (\pi(\phi_{M_0},\tau_{M_0})^w \otimes w(\nu_0)))\\
     &=\Ind_{M_1}^G (\Ind_{M_0'}^{M_1} (\pi({}^{\widehat{w}}\phi_{M_0},{}^{\widehat{w}}\tau_{M_0}) )\otimes w(\nu_0)).
\end{align*}
For the last equality, we use Desideratum \ref{desideratum: isomorphism} and the notation discussed below it. Now we decompose the tempered induction $\Ind_{M_0'}^{M_1}\pi({}^{\widehat{w}}\phi_{M_0},{}^{\widehat{w}}\tau_{M_0})$ using Desideratum \ref{desideratum: induction of tempered representations} explicitly. In the following, if $Q_1 \leq Q_2 \leq G$ are standard Levi subgroups, we let $\xi_{Q_1}^{Q_2}$ denote the $L$-embedding from ${}^L Q_1$ to ${}^L Q_2$.    Set $\phi'={\xi_{M_0'}^{M_1}} \circ {}^{\widehat{w}} \phi_{M_0}$. We have
\[\Ind_{M_0'}^{M_1}\pi({}^{\widehat{w}}\phi_{M_0},{}^{\widehat{w}}\tau_{M_0})= \bigoplus_{\tau' \in \widehat{A}_{\phi'} } m( {}^{\widehat{w}}\tau_{M_0},\tau') \pi(\phi', \tau').\]
Each $\pi(\phi',\tau')$ is a tempered representation of $M_1$ and $w(\nu_0)$ is positive with respect to $P_1$, so  $\Ind_{P_1}^G(\pi(\phi',\tau')\otimes w(\nu_0))$ is a standard representation. Put
\[\phi'':= \xi_{M_1}^G \circ \big(\phi' \otimes w(\nu_0)\big) = \xi_{M_0'}^G \circ \big({}^{\widehat{w}} \phi_{M_0}\otimes w(\nu_0)\big).\]
By Desideratum \ref{desideratum: Langlands classification} and Proposition \ref{prop: SZ18}(3) that standard module is $\Delta(\phi'',\tau'')$, where $\tau''\in\irrepdual{A}_{\phi''}$ corresponds to $\tau'$ under the isomorphism $A_{\phi'}\simeq A_{\phi''}$ induced by $\xi_{M_1}^G$. As $\xi_{M_0'}^G=\xi_{M_1}^G\circ \xi_{M_0'}^{M_1}$, this isomorphism is compatible with the two maps from $A_{{}^{\widehat{w}}\phi_{M_0}}$, so that $m({}^{\widehat{w}}\tau_{M_0},\tau')=m({}^{\widehat{w}}\tau_{M_0},\tau'')$. We obtain in the Grothendieck group that
\begin{align*}
    \Ind_M^G \Delta(\phi_M, \tau_M)&= \sum_{\tau' \in \widehat{A}_{\phi'} }   m( {}^{\widehat{w}}\tau_{M_0},\tau')  \Ind_{{M_1}}^G \big( \pi(\phi', \tau') \otimes w(\nu_0)\big)\\
    &=\sum_{\tau'' \in \widehat{A}_{\phi''} }   m( {}^{\widehat{w}}\tau_{M_0},\tau'')  \Delta(\phi'', \tau'').
\end{align*}

It remains to compare $\phi''$ with $\phi$. The Langlands classification of $\phi_M$ reads $\phi_M=\xi_{M_0}^M\circ(\phi_{M_0}\otimes\nu_0)$, so that $\phi=\xi_M^G\circ \phi_M=\xi_{M_0}^G\circ(\phi_{M_0}\otimes\nu_0)$. Comparing with the second expression for $\phi''$ above, we see that $\phi''={}^{\widehat{w}}\phi$. Hence $\phi''=\phi$ as $L$-parameters of $G$, and the isomorphism $A_{\phi''}\to A_\phi$ is conjugation by $\widehat{w}^{-1}$. Thus, we may rewrite
\begin{align*}
      \Ind_M^G \Delta(\phi_M, \tau_M)
    &=\sum_{\tau'' \in \widehat{A}_{\phi''} }   m( {}^{\widehat{w}}\tau_{M_0},\tau'')  \Delta(\phi'', \tau'')\\
    &=\sum_{\tau \in \widehat{A}_{\phi} }   m( {}^{\widehat{w}}\tau_{M_0},{}^{\widehat{w}}\tau)  \Delta(\phi, \tau)\\
    &=\sum_{\tau \in \widehat{A}_{\phi} }   m( \tau_{M_0},\tau)  \Delta(\phi, \tau)\\
    &=\sum_{\tau \in \widehat{A}_{\phi} }   m( \tau_{M},\tau)  \Delta(\phi, \tau).
\end{align*}
Here, the fourth equality comes from the isomorphism $A_{\phi_{M_0}}\simeq A_{\phi_M}$ by Proposition \ref{prop: SZ18}(3), which identifies $\tau_{M_0} \in \widehat{A}_{\phi_{M_0}}$ with $\tau_M \in \widehat{A}_{\phi_M}$. This proves \eqref{induction formula}.

{Finally we deduce the refinement. If $\kappa_M(\tau_M)=\gamma$ and $m(\tau_M,\tau)\neq 0$, then $\kappa_G(\tau)=H_\xi(\gamma)$ by Lemma \ref{lemma: relating galois colomology class of irreps of component groups}; so only such $\tau$ contribute to \eqref{induction formula}, which is \eqref{equation: induction formula form by form}, and the map $\Ind_{M_\gamma}^{G_\delta}$ is well defined. This completes the proof.}
\end{proof}

\begin{remark}
Though \eqref{induction formula} is an equation in $K\Pi_{\lambda}^{\pure}(G)$, the representations involved are in the same pure inner form of $G$. More precisely, let $\gamma:= \kappa_M(\tau_M)\in H^1(F,M)$ and let $\delta:= H_{\xi}(\gamma) \in H^1(F,G)$. By Lemma \ref{lemma: relating galois colomology class of irreps of component groups}, if $m(\tau_M, \tau) \neq 0$, then $\kappa_G(\tau)= \delta$. In conclusion, we may rewrite \eqref{induction formula} form-by-form:
    \begin{equation}\label{equation: induction formula form by form}
        \Ind_{M_\gamma}^{G_\delta}(\Delta(\phi_M,\tau_M))=\sum_{\tau \in \irrepdual{A}_\phi,~\kappa_G(\tau)=H_\xi(\gamma)}{m(\tau_M, \tau)} \Delta(\phi, \tau).
    \end{equation}
\end{remark}

\section{Parabolic induction as the transpose of geometric restriction}
Throughout this section, we are in the setting of \S \ref{subsection: setting}. In Section \ref{subsection: pure parabolic induction and restriction}  we show that $\Ind_M^G=\prescript{t}{}{\xi^*}$, where the transpose is taken with respect to the pairing \eqref{eqn:pairing} determined by $\LLC^{\mathfrak{e}}$. For $\gamma\in H^1(F,M)$ and $\delta=H_\xi(\gamma)$, we also offer the obvious refinement of this identity in Section \ref{subsection: parabolic induction form-by-form and lift}, by showing that $\Ind_{M_\gamma}^{G_\delta}$ is $\Lift_{M_\gamma}^{G_\delta}$, where the latter is defined in Section \ref{subsection: lift}.

\subsection{Pure parabolic induction and restriction}\label{subsection: pure parabolic induction and restriction}


\begin{prop}\label{proposition: pure parabolic induction and restriction}
    Let $G$ be a quasisplit connected reductive algebraic group with Levi subgroup $M$. We retain the notation from Section \ref{subsection: setting}. Assume that there exists an $\LLC^{\mathfrak{e}}$ for Levi subgroups of $G$ satisfying Desiderata \ref{desideratum: Langlands classification}, \ref{desideratum: induction of tempered representations}, \ref{desideratum: isomorphism}, and \ref{desideratum: Haines}. Further assume the $p$-adic Kazhdan-Lusztig Hypothesis \ref{hypothesis: KLH} for both $G$ and $M$. For $\pi\in K\Pi_{\lambda_M}^\pure(M)$ and $\mathcal{F}\in K\Per_{H_\lambda}(V_\lambda)$, we have
    \begin{equation}\label{equation: ind res adjunction}\langle \Ind_M^G(\pi),\mathcal{F}\rangle = \langle \pi,\xi^*(\mathcal{F})\rangle {.}\end{equation}
{This is the implication \textup{(a)}$\Rightarrow$\textup{(b)} of Theorem \ref{theorem: main theorem in the introduction}.}
\end{prop}
\begin{proof}
    It suffices to check the equality on the bases of standard representations of $K\Pi_{\lambda_M}^\pure(M)$ and of standard constructible sheaves of $K\Per_{H_\lambda}(V_\lambda)$, as provided by \eqref{equation: gr grp perverse sheaves} and \eqref{equation: gr grp of reps}. Consider the standard representation $\Delta(\phi_M,\tau_M)$ attached to the enhanced $L$-parameter $(\phi_M,\tau_M)$ of $M$, where $\tau_M\in \irrepdual{A}_{\phi_{M}}$. Let $\mathcal{F}=\stsh{\tau}{C}$, where $C$ is an $H_\lambda$-orbit in $V_\lambda$ and $\tau \in \irrepdual{A}_C$. From Theorem \ref{theorem: induction formula}, we have that 
    \[\Ind_{M}^{G} (\Delta(\phi_M,\tau_M))=\sum_{\tau' \in \widehat{A}_{\phi}} m(\tau_M, \tau') \Delta(\phi, \tau').\]
    This, combined with the identification $A_\phi=A_{C_\phi}$ from \eqref{eqn: equivalence of component group and eq fund group}, where $C_\phi$ is the $H_\lambda$-orbit in $V_\lambda$ corresponding to $\phi$, we get  
    \begin{align*}
    {\langle\Ind_M^G(\Delta(\phi_M,\tau_M)),\stsh{\tau}{C}\rangle}
        &=   \displaystyle \sum_{\tau' \in \widehat{A}_{C_\phi}} m(\tau_M, \tau')\langle \Delta(\phi, \tau'), \stsh{\tau}{C}\rangle \\
          &=\begin{cases}
             m(\tau_M,\tau)\langle \Delta(\phi,\tau),\stsh{\tau}{C} \rangle & \text{ if } C=C_\phi,\\
             0 & \text{otherwise}.
          \end{cases}  \\
          &=\begin{cases}
             m(\tau_M,\tau) e(\kappa_G(\tau)) & \text{ if } C=C_\phi,\\
             0 & \text{otherwise}. 
          \end{cases}  
    \end{align*}
where the last equality is evaluated using the Kazhdan-Lusztig Hypothesis (Hypothesis \ref{hypothesis: KLH}). This computes the left hand side of \eqref{equation: ind res adjunction}.

Now we compute the right hand side of \eqref{equation: ind res adjunction}.
From  \eqref{eqn: restriction formula on standards}, we have
\[\xi^*(\stsh{\tau}{C}) = \sum_{\xi(D)\subset C,~\sigma \in \irrepdual{A}_D} m(\sigma,\tau) \cdot \stsh{\sigma}{D}.\] 
We obtain
\begin{align*}
     \langle \Delta(\phi_M,\tau_M), {\xi^*(\stsh{\tau}{C})}\rangle &= \displaystyle\sum_{\xi(D)\subset C,~\sigma \in \irrepdual{A}_D}m(\sigma,\tau)\langle \Delta(\phi_M,\tau_M),\stsh{\sigma}{D}\rangle \\&=\begin{cases}
         m(\tau_M,\tau) e(\kappa_M(\tau_M)) & \text{ if } C=C_\phi,\\
         0 & \text{ otherwise.}
         \end{cases}
\end{align*}
Here the second equality holds because a summand $\langle \Delta(\phi_M,\tau_M),\stsh{\sigma}{D}\rangle$ is nonzero only for $(D,\sigma)=(C_{\phi_M},\tau_M)$, which forces $C=C_\phi$, the saturation of the orbit $C_{\phi_M}$ in $V_\lambda$. The resulting value is then read off once again using Hypothesis \ref{hypothesis: KLH}. If $m(\tau_M,\tau)\neq 0$, then $\kappa_G(\tau)=H_\xi(\kappa_M(\tau_M))$ by Lemma \ref{lemma: relating galois colomology class of irreps of component groups}, whence $e(\kappa_G(\tau))=e(\kappa_M(\tau_M))$ by \eqref{eqn Kottwitz sign Levi}. 
\end{proof}

\subsection{Parabolic induction form-by-form and lift}\label{subsection: parabolic induction form-by-form and lift}

Recall the map \[\Lift_{M_\gamma}^{G_\delta}: K\Pi_{\lambda_M}(M_\gamma)\to K\Pi_\lambda(G_\delta)\] from Section \ref{subsection: lift}.
\begin{prop}\label{lemma: lift is induction}
    Assume that there exists a $\LLC^{\mathfrak{e}}$ for Levi subgroups of $G$ satisfying Desiderata \ref{desideratum: Langlands classification}, \ref{desideratum: induction of tempered representations}, \ref{desideratum: isomorphism}, and \ref{desideratum: Haines}. Further assume the $p$-adic Kazhdan-Lusztig Hypothesis \ref{hypothesis: KLH} both for $G$ and $M$. We retain the notation from Sections \ref{subsection: setting} and \ref{subsection: lift}. We have
    \[ \Ind_{M_\gamma}^{G_\delta} = \Lift_{M_\gamma}^{G_\delta}.\]
\end{prop}
\begin{proof}
    It suffices to show that 
    \begin{equation}\label{eqn: need to show}
        \langle \Ind_{M_\gamma}^{G_\delta}(\Delta(\phi_M, \tau_M)), \stsh{\tau}{C} \rangle = \langle \Lift_{M_\gamma}^{G_\delta}(\Delta(\phi_M, \tau_M)), \stsh{\tau}{C}\rangle
    \end{equation}
    for all $ (\phi_M,\tau_M) \in \Phi^{\mfr{e}}(M)$ such that $\kappa_M(\tau_M)=\gamma$, and all orbits $C \in V_{\lambda}/H_{\lambda}$ and $\tau \in \irrepdual{A}_{C}$.

    For the left hand side of \eqref{eqn: need to show}, using  \eqref{equation: induction formula form by form}, we have 
    \begin{align*}
        \langle \Ind_{M_\gamma}^{G_\delta}(\Delta(\phi_M, \tau_M)),\stsh{\tau}{C} \rangle &= \sum_{\tau' \in \irrepdual{A}_\phi,~\kappa_G(\tau')=\delta} m(\tau_M,\tau')\langle \Delta(\phi, \tau'), \stsh{\tau}{C} \rangle \\
        &= \begin{cases}
            m(\tau_M,\tau)\langle \Delta(\phi, {\tau}), \stsh{\tau}{C_\phi} \rangle  & \text{if } C=C_\phi \text{ and } \kappa_G(\tau)=\delta,\\
            0 & \text{otherwise}.
        \end{cases}\\
        &= \begin{cases}
            m(\tau_M,\tau)e_G(\delta)  & \text{if } C=C_\phi \text{ and } \kappa_G(\tau)=\delta,\\
            0 & \text{otherwise}.
        \end{cases}
    \end{align*}
The computation of $\langle \Delta(\phi, {\tau}), \stsh{\tau}{C_\phi} \rangle$ in {the} last equality follows from the Kazhdan-Lusztig Hypothesis \ref{hypothesis: KLH}. 

For the right hand side of \eqref{eqn: need to show}, using \eqref{equation: lift}, we have
\begin{align*}
   \langle \Lift_{M_\gamma}^{G_\delta}&(\Delta(\phi_M, \tau_M)), \stsh{\tau}{C}\rangle\\
   &= \big\langle \pr^G_\delta \circ \prescript{t}{}{\xi^*} \circ \iota^M_\gamma (\Delta(\phi_M, \tau_M)), \stsh{\tau}{C}\big\rangle\\[2pt]
   &= \big\langle \pr^G_\delta \circ \prescript{t}{}{\xi^*} (\Delta(\phi_M,\tau_M)),\stsh{\tau}{C}\big\rangle
      && \text{by }\eqref{equation: inclusion of grothendieck groups of reps}\\[2pt]
   &= \begin{cases}
        \big\langle \prescript{t}{}{\xi^*} (\Delta(\phi_M,\tau_M)), \stsh{\tau}{C}\big\rangle & \text{if }\kappa(\tau)=\delta,\\
        0 & \text{otherwise},
      \end{cases}
      && \text{by }\eqref{equation: transpose of projection}\\[2pt]
   &= \begin{cases}
        \big\langle \Delta(\phi_M,\tau_M), \xi^*(\stsh{\tau}{C})\big\rangle & \text{if }\kappa(\tau)=\delta,\\
        0 & \text{otherwise},
      \end{cases}
      && \text{by }\eqref{equation: transpose of restriction}\\[2pt]
   &= \begin{cases}
        m(\tau_M,\tau)\, e_M(\gamma) & \text{if } C=C_\phi \text{ and } \kappa_G(\tau)=\delta,\\
        0 & \text{otherwise},
      \end{cases}
      && \text{as in Proposition \ref{proposition: pure parabolic induction and restriction}.}
\end{align*}
Finally, note that $e_M(\gamma)= e_{G}(H_{\xi}(\gamma))= e_G(\delta)$ by \eqref{eqn Kottwitz sign Levi}. We conclude that \eqref{eqn: need to show} holds. 
\end{proof}

\section{The converse implication: recovering the two conjectures}\label{section: converse}
In this section, we prove the implication \textup{(b)}$\Rightarrow$\textup{(a)} of Theorem \ref{theorem: main theorem in the introduction}. It is the converse of Proposition \ref{proposition: pure parabolic induction and restriction} in the following sense. Recall the maps $\Ind_M^G:K\Pi_{\lambda_M}^\pure(M)\to K\Pi_\lambda^\pure(G)$ from \eqref{equation: induction definition} and $\xi^*:K\Per_{H_\lambda}(V_\lambda)\to K\Per_{H_{\lambda_M}}(V_{\lambda_M})$ from \eqref{equation: restriction linear map perverse sheaves}. Assuming the identity $\Ind_M^G= {}^t \xi^{\ast}$,  we recover both the $p$-adic analogue of the Kazhdan-Lusztig hypothesis and the decomposition of tempered induction in
Desideratum \ref{desideratum: induction of tempered representations}. 

\begin{prop}\label{theorem: ind res adjunction implies klh}
Let $G$ be a quasisplit connected reductive algebraic group.  Assume that there exists an $\LLC^{\mathfrak{e}}$ for Levi subgroups of $G$ satisfying Desiderata \ref{desideratum: Langlands classification} and \ref{desideratum: Haines}. Suppose further that for any standard Levi subgroup $M$ and infinitesimal parameter $\lambda_M$, the identity \eqref{equation: ind res adjunction}
\[\langle \Ind_M^G(\pi),\mathcal{F}\rangle=\langle \pi, \xi^*(\mathcal{F})\rangle  \]
holds for all $\pi \in K\Pi_{\lambda_M}^{\pure}(M)$ and $\mathcal{F} \in K \Per_{H_{\lambda}}(V_{\lambda
})$. Then
\begin{enumerate}
    \item [(a)] The $p$-adic Kazhdan-Lusztig hypothesis holds$:$ 
\begin{equation*}
    \begin{array}{rcl}
        K\Pi_\lambda^\pure(G) \times K\Per_{H_\lambda}(V_\lambda)
            &\to& \mathbb{Z}\\
        \langle \Delta(\pi), \tau_{C_{\phi}}^{} \rangle 
            &\mapsto& 
            \begin{cases}
                e(\pi) & \LLC^{\mathfrak{e}}_G(\pi) = (\phi,\tau), \\
                0 & \text{else}.
            \end{cases}
    \end{array}
\end{equation*}
    \item [(b)]  Desideratum \ref{desideratum: induction of tempered representations} holds: for any standard Levi subgroup $M$ of $G$, and tempered $L$-parameter $\phi_M$ of $M$,
    \[\Ind_{M}^{G} \big(\pi(\phi_M,\tau_M)\big)=\bigoplus_{\tau \in \irrepdual{A}_{\phi}} m(\tau_M, \tau)\, \pi(\phi, \tau).\]
\end{enumerate}
\end{prop}
\begin{proof}
We first prove Part (a).   For an irreducible representation $\pi \in \Pi_{\lambda}^\pure(G)$, write 
    \[\Delta(\pi)=\Ind_M^G(\pi_M) \in K\Pi_{\lambda}^\pure(G),\]
    where $\pi_M$ is an essentially tempered representation in $\Pi_{\lambda_M}^\pure(M)$ with enhanced Langlands parameter $(\phi_M, \sigma)$.  Let $\LLC_G^{\mathfrak{e}}(\pi)=(\phi,\tau)$; then $\phi=\xi\circ\phi_M$ as in \eqref{equation: setting L par factors}, and in particular $\phi$ has infinitesimal parameter $\lambda= \xi \circ \lambda_M$. Moreover, we know from Proposition \ref{prop: SZ18}(3) that 
\begin{equation}\label{equation: isom of comp groups tempered}
    A_\phi \simeq A_{\phi_M}, 
\end{equation}
so we may write \begin{equation}\label{equation: sigma is tau}
  \tau \circ   A_{\phi_M}^{\phi}=\sigma, \text{ or with mild abuse of notation, } \tau=\sigma.
\end{equation}
Let $C_{\phi_M}$ denote the $H_{\lambda_M}$-orbit in $V_{\lambda_M}$ corresponding to $\phi_M$ according to \eqref{eqn: vogan langlands correspondence}. The orbit $C_{\phi_M}$ is open in $V_{\lambda_M}$; see Remark \ref{rmk LC tempered}. Let $C$ be an $H_\lambda$-orbit in $V_\lambda$ and $\gamma \in \irrepdual{A}_C$. By assumption of identity \eqref{equation: ind res adjunction},
    \begin{equation}\label{equation: application of ind res adj}
        \langle \Ind_M^G(\pi_M), \stsh{\gamma}{C}\rangle = \langle \pi_M, \xi^*(\stsh{\gamma}{C})\rangle. \end{equation}
By \eqref{eqn: restriction formula on standards}, the sheaf $\xi^*(\stsh{\gamma}{C})$ is supported on those $H_{\lambda_M}$-orbits $D$ with $d\xi(D)\subseteq C$. Since the pairing on the right-hand side of \eqref{equation: application of ind res adj} only sees the orbit $C_{\phi_M}$ attached to $\pi_M$, it vanishes unless $d\xi(C_{\phi_M})\subseteq C$. 
The orbit $C_\phi$ is by definition the saturation of $C_{\phi_M}$ in $V_\lambda$ so $C=C_\phi$.
Recall from \eqref{eqn: equivalence of component group and eq fund group} that $A_{\phi}$ (resp. $A_{\phi_M}$) is the equivariant fundamental group of $C_\phi$ (resp. $C_{\phi_M}$). Using the formula \eqref{eqn: restriction formula on standards} for $\xi^*$ with the identification \eqref{equation: isom of comp groups tempered} in mind, we get
\[\langle \pi_M, \xi^*(\stsh{\gamma}{C})\rangle= \langle \pi_M,\stsh{\gamma}{C_{\phi_M}} \rangle.\]
Let $j:C_{\phi_M}\hookrightarrow V_{\lambda_M}$ denote the inclusion. Since $C_{\phi_M}$ is open in $V_{\lambda_M}$ ($\ol{C}_{\phi_M}=V_{\lambda_M}$), and $\locsys{\gamma}{C_{\phi_M}}[\dim C_{\phi_M}]\in \Per(C_{\phi_M})$, we have from \eqref{equation: standard sheaf as a sum of perverse sheaf and boundary terms} 
\[(-1)^{\dim C_{\phi_M}}\stsh{\gamma}{C_{\phi_M}}=\mathcal{IC}(\locsys{\gamma}{C_{\phi_M}})+ \text{ boundary terms} ~~\in K\Per_{H_{\lambda_M}}(V_{\lambda_M}).\]

Substituting $\stsh{\gamma}{C_{\phi_M}}$ in the pairing, the boundary terms pair trivially with $\pi_M$, so we get

\begin{align}\label{eq pairing converse}
\nonumber \langle \pi_M, \xi^*(\stsh{\gamma}{C_\phi})\rangle 
&= \left\langle \pi_M, \stsh{\gamma}{C_{\phi_M}}\right\rangle \\
&= (-1)^{\dim C_{\phi_M}}\left\langle \pi_M, \mathcal{IC}\left(\locsys{\gamma}{C_{\phi_M}}\right)\right\rangle \\
 \nonumber &=\begin{cases}
    (-1)^{2\dim C_{\phi_M}}~e(\pi_M) & \text{ if } C=C_\phi \text{ and } \gamma=\sigma,\\
    0 & \text{otherwise,} 
\end{cases}\\
\nonumber   &= \begin{cases}
    e(\pi) & \text{ if } C=C_\phi \text{ and } \gamma=\tau,\\
    0 & \text{otherwise.}
    \end{cases}
\end{align}
where the second to last equality follows from the definition of the pairing \eqref{eqn:pairing}, and in the last equality the condition follows from \eqref{equation: sigma is tau} and $e(\pi_M)=e(\pi)$  from \eqref{eqn Kottwitz sign Levi}.

In summary, we have computed that
\[ \langle \Delta(\pi), \stsh{\gamma}{C} \rangle=\begin{cases}
    e(\pi) & \text{ if } C=C_\phi \text{ and } \gamma=\tau,\\
    0 & \text{otherwise.}
    \end{cases} \]
In the nonzero instance, the condition $\langle \Delta(\pi),{\stsh{\gamma}{C}}\rangle = e(\pi) $ if and only if $ (C,\gamma)=(C_\phi,\tau) $ can be rewritten as $\langle \Delta(\pi),\stsh{\tau}{C_\phi}\rangle = e(\pi) $ if and only if $ \LLC_G^{\mathfrak{e}}(\pi)=(\phi, \tau)$. This completes the proof of Part (a).

Now we prove Part (b). From now on, $\phi_M$ is a tempered $L$-parameter of $M$ with infinitesimal parameter $\lambda_M$, and $\phi=\xi \circ \phi_M, \lambda= \xi \circ \lambda_M$ are the corresponding tempered $L$-parameter and infinitesimal parameter of $G$. It is well known that induction of a tempered representation is a semisimple direct sum of tempered representations. Moreover, by Desideratum \ref{desideratum: Haines}, we know that 
\[ \Ind_{M}^G (\pi(\phi_M,\tau_M))= \bigoplus_{\tau \in \irrepdual{A}_{\phi}} m_{\tau} \pi(\phi,\tau),\]
where $m_{\tau}$ are non-negative integers. Indeed, if $\pi=\pi(\phi',\tau')$ appears in $\Ind_M^G(\pi(\phi_M,\tau_M))$, then 
\[ \lambda_{\phi'}= \xi \circ \lambda_{\phi_M}= \lambda_{\phi}=:\lambda.\]
Since $\phi'$ and $\phi$ are both tempered, they correspond to the unique open orbit in $V_{\lambda}$ (\cite{cunningham2024whittakernormalizationpadicabvpackets}*{Proposition 3.2}). Hence, $\phi'$ is equal to $\phi$. 

It thus remains to show that the multiplicity $m_{\tau}$ matches $ m(\tau_M,\tau)$. Now we have 
\begin{align*}
    m_{\tau}& =  e(\tau) (-1)^{\dim(C_{\phi})} \langle \Ind_M^G (\pi(\phi_M,\tau_M)) , \IC(\locsys{\tau}{C_{\phi}}) \rangle\\
    &=e(\tau) \langle \Ind_M^G (\pi(\phi_M,\tau_M)) , \stsh{\tau}{C_{\phi}} \rangle \\
    &=e(\tau)  \langle \pi(\phi_{M}, \tau_M) , \xi^{\ast} (\stsh{\tau}{C_{\phi}})\rangle\\
    &= e(\tau)  \left\langle \pi(\phi_{M}, \tau_M) , \sum_{\xi (D) \subseteq C_{\phi},\ \sigma \in \irrepdual{A}_{D}} m(\sigma, \tau) \stsh{\sigma}{D}\right\rangle\\
    &= e(\tau) e(\tau_M) m(\tau_M, \tau)\\
    &= m(\tau_M, \tau).
\end{align*}
Here, the second and the fifth equalities follow from the same reason as \eqref{eq pairing converse}; the third equality is our assumption \eqref{equation: ind res adjunction}; the fourth equality is \eqref{eqn: restriction formula on standards}; and the last equality is \eqref{eqn Kottwitz sign Levi}. This completes the proof of Part (b) and the proposition.
\end{proof}

\section{Supercuspidal sheaves}\label{section: supercuspidal sheaves}
In this section, we use the geometric description of parabolic induction \eqref{equation: ind res adjunction} to 
describe the simple perverse sheaves corresponding to supercuspidal representations under local Langlands correspondence \eqref{eqn: vogan langlands correspondence}.

\begin{defn}\label{def sc sheaf}
A simple perverse sheaf $\mathcal{P} \in \Per_{H_\lambda}(V_\lambda)$ is called  \emph{supercuspidal} if for every proper Levi subgroup $M$ for which $\lambda$ factors through the embedding $\xi: {}^L M \hookrightarrow {}^L G$ as $\lambda_M$, we have
\begin{equation}\label{equation: supercuspidal sheaf condition}
    \xi^*(\mathcal{P})=0 \quad \text{in } K\Per_{H_{\lambda_M}}(V_{\lambda_M}).
\end{equation}
\end{defn}

\begin{prop}
    Keep the assumptions in Proposition \ref{proposition: pure parabolic induction and restriction} for $G$. A representation $\pi \in \Pi_{\lambda}^{\pure}(G)$ is supercuspidal if and only if the corresponding simple perverse sheaf $\mathcal{P}(\pi)\in \Per_{H_\lambda}(V_\lambda)$ is supercuspidal.
\end{prop}
\begin{proof}
    We prove that $\pi$ is not supercuspidal if and only if $\mathcal{P}(\pi)$ is not supercuspidal.

    Suppose that $\pi$ is not supercuspidal, then there exists a proper Levi subgroup $M$ and $\sigma \in \Pi^{\pure}(M)$ such that $\pi$ is a subquotient of $\Ind_M^G(\sigma)$.  Namely, 
    \[0 \neq \langle \Ind_M^G(\sigma), \mathcal{P}(\pi) \rangle.\]
    Write $\lambda_M$ for the infinitesimal parameter of $\sigma$ and write $\xi: {}^LM \hookrightarrow {}^L G$ for the embedding. Then Desideratum \ref{desideratum: Haines} implies that $({}^L M \hookrightarrow {}^L G) \circ \lambda_M= \lambda$. Now Proposition \ref{proposition: pure parabolic induction and restriction} gives
    \[  0 \neq \langle \Ind_M^G(\sigma), \mathcal{P}(\pi) \rangle=  \langle \sigma, \xi^{\ast}(\mathcal{P}(\pi)) \rangle. \]
    Therefore, we must have $ \xi^{\ast}(\mathcal{P}(\pi))\neq 0$ and $\mathcal{P}(\pi)$ is not supercuspidal.

    Conversely, suppose that $\mathcal{P}(\pi)$ is not supercuspidal. Take a Levi subgroup $M$ such that $\xi^{\ast} (\mathcal{P}(\pi)) \neq 0$ where $\xi:{}^L M \hookrightarrow {}^L G$. Then there exists a $\sigma \in \Pi_{\lambda_M}^{\pure}(M)$ such that $\langle \sigma ,\xi^{\ast} (\mathcal{P}(\pi))  \rangle \neq 0 $. Then again by Proposition \ref{proposition: pure parabolic induction and restriction}, we obtain
    \[ 0 \neq \langle \sigma ,\xi^{\ast} (\mathcal{P}(\pi))\rangle = \langle \Ind_M^G \sigma , \mathcal{P}(\pi) \rangle. \]
    This implies that $\pi$ appears in the induction $\Ind_M^G \sigma$, and hence $\pi$ is not supercuspidal. This completes the proof of the proposition.
\end{proof}

\begin{remark} If $\lambda$ does not factor through the $L$-group of any proper Levi subgroups of $G$, then every simple perverse sheaf in $\Per_{H_\lambda}(V_{\lambda})^{\mathrm{simple}}_{\mathrm{/iso}}$ is supercuspidal by Definition \ref{def sc sheaf}. In this case, all of the representations in $\Pi_{\lambda}^{\pure}(G)$ (which turns out to be an $L$-packet with $L$-parameter $\lambda$) are supercuspidal by Desideratum \ref{desideratum: Haines} as in the proof above.   
\end{remark}

\bibliographystyle{amsalpha} 
\bibliography{references} 
\end{document}